\documentclass{amsart}
\usepackage{amsfonts,  bbm, array}
\usepackage{color}

\usepackage{overpic}

\newtheorem{theorem}{Theorem}[section]
\theoremstyle{plain}

\newtheorem{conj}[theorem]{Conjecture}
\theoremstyle{remark}

\numberwithin{equation}{section}

\newcommand{\re}{\operatorname{Re}}
\newcommand{\im}{\operatorname{Im}}
\newcommand{\tr}{\operatorname{tr}}

\newcommand{\vol}{\operatorname{vol}}

\newcommand{\res}{\operatorname{Res}}
\newcommand{\rank}{\operatorname{rank}}

\newcommand{\brak}[1]{\langle #1 \rangle}
\newcommand{\abs}[1]{\left|#1 \right|}

\newcommand{\bbR}{\mathbb{R}}
\newcommand{\bbH}{\mathbb{H}}
\newcommand{\bbC}{\mathbb{C}}
\newcommand{\bbZ}{\mathbb{Z}}
\newcommand{\bbN}{\mathbb{N}}

\newcommand{\calR}{\mathcal{R}}
\newcommand{\calH}{\mathcal{H}}

\newcommand{\calL}{\mathcal{L}}

\newcommand{\calW}{\mathcal{W}}

\newcommand{\vep}{\varepsilon}

\begin{document}

\title{Distribution of resonances for hyperbolic surfaces} 
\author[D. Borthwick]{David Borthwick}
\address{Department of Mathematics and Computer Science, Emory
University, Atlanta, Georgia, 30322, USA}
\thanks{Supported in part by NSF\ grant DMS-0901937.}
\email{davidb@mathcs.emory.edu}
\date{May 19, 2013}

\begin{abstract}
We study the distribution of resonances for geometrically finite hyperbolic surfaces of infinite area by counting resonances numerically.  The resonances are computed as zeros of the Selberg zeta function, using an algorithm for computation of the zeta function for Schottky groups.   Our particular focus is on three aspects of the resonance distribution that have attracted attention recently: the fractal Weyl law, the spectral gap, and the concentration of decay  rates.  
\end{abstract}

\maketitle
\tableofcontents
\section{Introduction}

A smooth, geometrically finite hyperbolic surface $X = \Gamma\backslash\bbH$ has finite genus, with ends consist of a finite number of hyperbolic funnels and/or cusps.  We'll assume that $X$ has infinite area, so there is at least one funnel.  Under this assumption, the (positive) Laplacian $\Delta_X$ has absolutely continuous spectrum $[\tfrac14,\infty)$ and finitely many discrete eigenvalues in $(0,\tfrac14)$, with no embedded eigenvalues.  The resolvent $R_X(s) := (\Delta_X - s(1-s))^{-1}$ is well-defined for $\re s > \tfrac12$, as long as $s(1-s)$ is not an $L^2$ eigenvalue of $\Delta_X$.  By Mazzeo-Melrose \cite{MM:1987} and Guillop\'e-Zworski \cite{GZ:1995a}, $R_X(s)$ admits a meromorphic continuation to $s \in \bbC$, with poles of finite rank.  The \emph{resonances} of $X$ are defined as the poles of this continuation.  We let $\calR_X$ denote the set of resonances repeated according to the multiplicity given by the rank of the residue of $R_X(s)$ at the pole.

\begin{figure}
\begin{tabular}{cc}
\begin{overpic}[scale=.7]{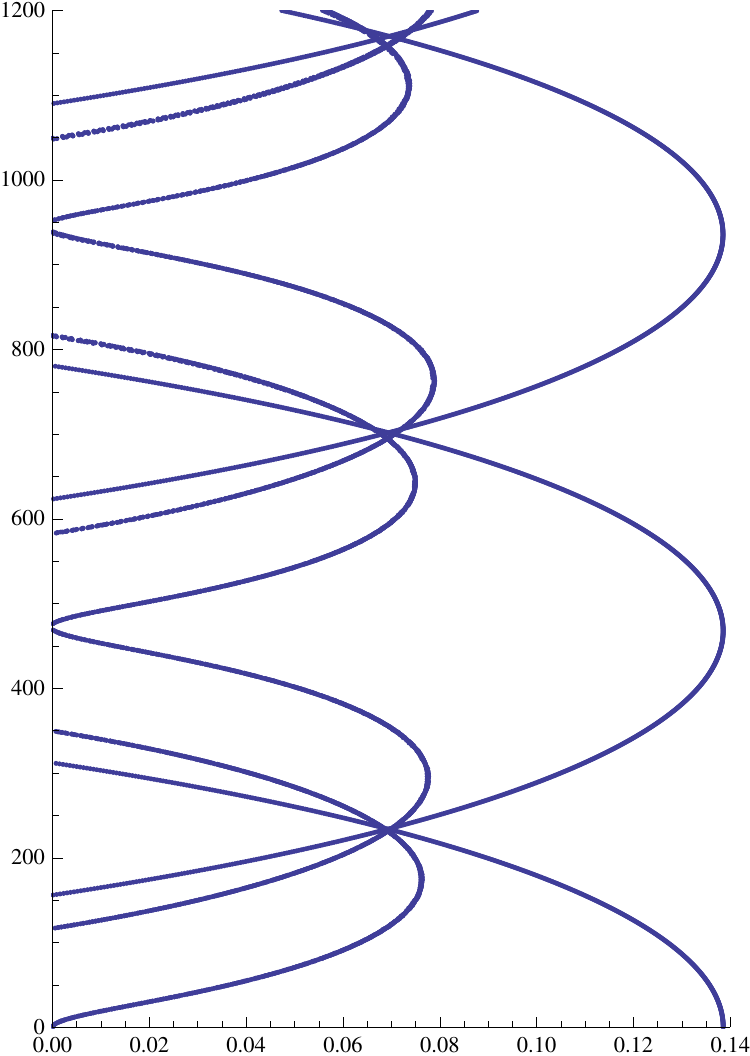}
\put(25,-3){$\scriptstyle X(10,10,10)$}
\end{overpic}  &
\begin{overpic}[scale=.7]{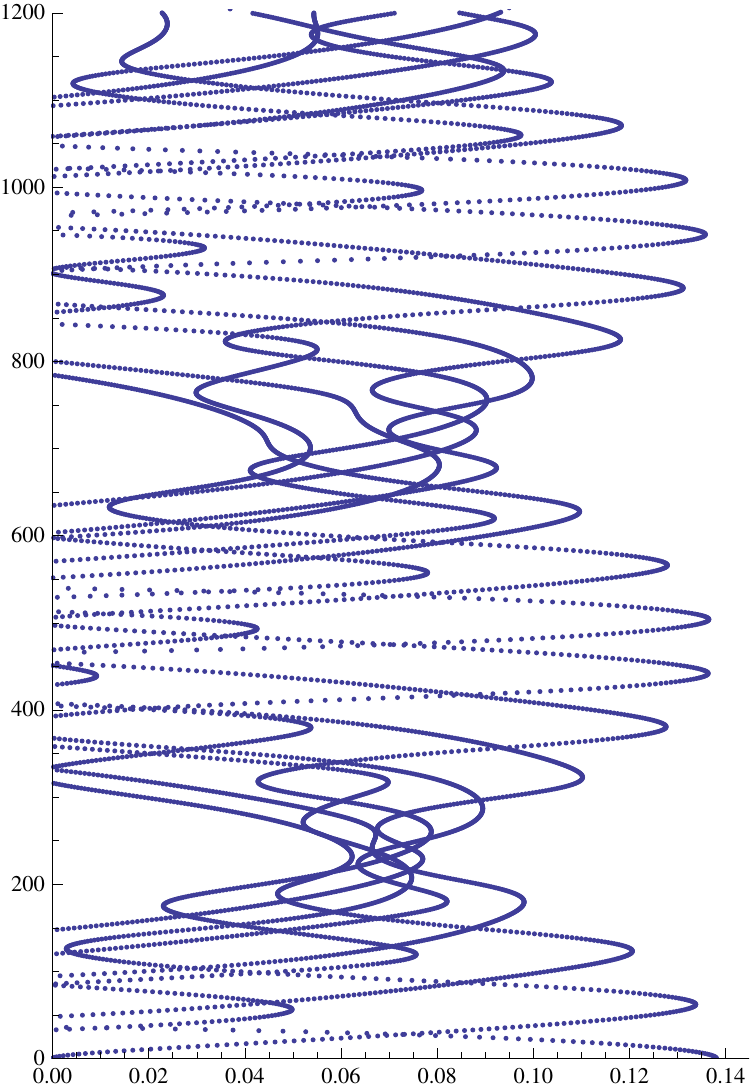}
\put(25,-3){$\scriptstyle X(10,10,10.1)$}
\end{overpic} 
\end{tabular}
\caption{Resonances of 3-funnel surfaces $X(\ell_1,\ell_2,\ell_3)$, consisting of funnels attached to a hyperbolic pair of pants with boundary lengths $\ell_1,\ell_2,\ell_3$.}\label{ResPlot101.fig}
\end{figure}

Associated to the surface (or the group) is the Selberg zeta function, $Z_X(s)$, defined as a product over the primitive length spectrum $\calL_X$, 
\begin{equation}\label{Z.def}
Z_X(s) := \prod_{\ell\in\calL_X} \prod_{k=0}^\infty \left( 1 - e^{(s+k)\ell} \right),
\end{equation}
converging uniformly on compact sets for $\re s \ge 1$.
The zeta function admits meromorphic continuation to $s\in\bbC$ for any geometrically finite hyperbolic surface.  This follows from the Selberg trace formula if the area is finite, and was extended to the general case by Guillop\'e \cite{Guillope:1992}.  By Patterson-Perry \cite{PP:2001} and Borthwick-Judge-Perry \cite{BJP:2005}, the divisor of  $Z_X(s)$ consists of zeros at the resonance set $\calR_X$, plus a set of ``topological'' zeros (depending on the Euler characteristic of $X$) and poles (depending on the number of cusps.)  In particular, $Z_X(s)$ is entire if and only if $X$ has no cusps (equivalent to $\Gamma$ being convex co-compact).
See \cite{Borthwick} for additional background on this link between the resonances and the zeta function.

Figure~\ref{ResPlot101.fig} shows sample plots of the resonance sets of two hyperbolic surfaces, computed as zeros of the corresponding zeta function by the method described in \S\ref{zeta.sec}--\ref{comp.sec}.  The curves in the plots are made up of discrete points;  the spacing between resonances is small enough that individual resonances are not necessarily resolved at this scale.  The plots show only $\re s\ge 0$ because the resonance sets are symmetric with respect to conjugation, by self-adjointness.
In our definition of resonance we use the parameter $s$ because of the connection to the Selberg zeta function.  A more traditional choice from the point of view of scattering theory would be $s = \tfrac12 + i\lambda$, which would flip the axes in the plots.

The theory of the zeta function for infinite-area surfaces is of course an extension of the much older Selberg theory for compact or finite-area hyperbolic surfaces.   In the finite-area case, Selberg \cite{Selberg:1989} showed that the resonances are confined to a vertical strip $-c_1 \le \re s \le 1$.   Furthermore, the Selberg trace formula leads a Weyl law for resonances of the form,
\begin{equation}\label{finite.weyl}
\#\left\{\zeta\in\calR_X:\> \abs{\im \zeta} \le t\right\} \sim  \frac{\vol(X)}{2\pi} t^2,
\end{equation}
proven in this case by M\"uller \cite{Muller:1992} and Parnovski \cite{Parnovski:1995}.
Another result of Selberg \cite{Selberg:1990} shows that the resonances of a finite-area hyperbolic surface accumulate on the critical line $\re s = \tfrac12$, in the sense that for any $\vep>0$, there are infinitely resonances in the strip $\frac12 -\vep < \re s < \frac12$.  Dyatlov \cite[eq.~(1.17)]{Dyatlov:2012} recently pointed out that under the assumption of finitely many embedded eigenvalues (a condition which is known to hold generically), most resonances lie near the critical line in the sense that
\[
\#\left\{\zeta\in\calR_X:\> \abs{\zeta - \tfrac12} \le t, \>\re s < \tfrac12 - \vep\right\} = o(t^2),
\]
for any $\vep>0$.  

By comparison, in the infinite-area case a great deal is still unknown about the distribution of resonances.  For example, there is no Weyl law analogous to \eqref{finite.weyl}.  Guillop\'e-Zworski \cite{GZ:1995a, GZ:1997} proved that for non-elementary hyperbolic surfaces (as well as compactly supported perturbations), 
\begin{equation}\label{GZ.bnd}
\#\left\{\zeta\in\calR_X:\> \abs{\zeta-\tfrac12} \le t\right\} \asymp t^2.
\end{equation}
In \cite{Borthwick:2012}, we showed that there exist constants for these upper and lower that depend only on $\chi(X)$ and the sum of lengths of the geodesic funnel boundaries, with the upper constant being sharp in the sense that it reduces to \eqref{finite.weyl} in the finite-area case. 

The distribution of resonances for infinite-area surfaces has attracted renewed attention lately on the arithmetic side, in the work of Bourgain-Gamburd-Sarnak \cite{BGS:2011}.   Infinite-area surfaces are not included in the classical Selberg theory, since arithmetic Fuchsian groups yield quotients of finite volume.  However, the main point in \cite{BGS:2011} was a generalization of Selberg's famous $\tfrac{3}{16}$ theorem to infinite-area ``congruence'' surfaces.   This amounts to a uniform estimate on the spectral gap between the first (rightmost) resonance and the bulk of $\calR_X$.  Such estimates have applications to upper bounds for the affine linear sieve.

Another source of recent interest in the distribution of resonances for infinite-area hyperbolic surfaces comes from the 
theory of chaotic scattering systems in physics, as exemplified by $n$-disk scattering systems.  In the last decade a number of interesting conjectures have arisen for such systems, involving connections between the distribution of resonances near the critical line ($\re s = \tfrac12$ in our case) and chaotic dynamics of the associated classical system.  The conjectures have been supported by both numerical and experimental evidence, and also by mathematical theorems in some cases.  

With these recent developments as motivation, we present here a numerical investigation 
in the spirit of Hejhal-Rackner \cite{HR:1992}.  We will focus on some particular aspects of the distribution of resonances for hyperbolic surfaces:
\begin{enumerate}
\item  The asymptotic density of resonances near the critical line: the fractal Weyl law.
\item  The spectral gap between the first resonance and the bulk of the resonance set.
\item  Concentration of imaginary parts (quantum decay rates) at half of the classical escape rate.
\end{enumerate}
We will give further details of the specific conjectures, as well as background and references, in the relevant sections \S\ref{frweyl.sec}--\ref{decay.sec} below.   For the finite-area case, the conjectures all correspond to known theorems, by the results mentioned above.

The dimension $\delta$ of the limit set of $\Gamma$ plays a central role in our investigations.  Patterson's famous result \cite{Patterson:1976} equated the dimension of the limit set with the exponent of convergence of $\Gamma$, and showed that $\delta(1-\delta)$ was the lowest eigenvalue of the Laplacian if $\delta>\tfrac12$.  Patterson extended this result to $\delta\le \tfrac12$ in \cite{Patterson:1988}, showing that $\delta$ is always the location of the first resonance.  In addition, the value of $\delta$ determines the dimension for the fractal Weyl law and the classical escape rate, so it is the crucial parameter in all three conjectures mentioned above.

The basis for our computations is an algorithm for computing the Selberg zeta function for a classical Schottky group from Jenkinson-Pollicott \cite{JP:2002}.  By a result of Button \cite{Button:1998}, any convex co-compact Fuchsian group is representable as a classical Schottky group, so in principle this technique would allow the approximation of  resonances for any smooth hyperbolic surface without cusps.  In practice, however, the number of computations required increases exponentially with the number of generators of the Schottky group.  To keep the computations manageable we have confined our attention to groups with 2 generators, which limits us to surfaces with $\chi = -1$.  Another limitation on the method is a severe degradation of the convergence rate for $\re s < 0$.  We are thus not able to investigate global asymptotics of the form \eqref{finite.weyl}.  Fortunately, the resonances of greatest interest are those nearest the critical line, and the convergence rate is stable up to very large values of $\im s$, provided we choose group parameters so that $\delta$ is relatively small.

\vskip12pt\noindent
\textbf{Acknowledgment}.  I would like to thank Maciej Zworski for encouragement to develop this project, and for helpful comments and suggestions along the way.

\section{Spectral theory preliminaries}\label{prelim.sec}

As noted in the introduction, for a geometrically finite hyperbolic surface $X = \Gamma\backslash \bbH$ the resolvent is written $R_X(s) := (\Delta - s(n-s))^{-1}$ and admits a meromorphic continuation to $s\in\bbC$.   The resonance set $\calR_X$ is then defined to be the set of poles of $R_X(s)$, with multiplicities given by
\[
m_\zeta := \rank \res_{s=\zeta} R_X(s).
\]
The corresponding global counting function is
\[
N_X(t) := \#\bigl\{\zeta\in\calR_X:\>\abs{\zeta-\tfrac12} \le t\bigr\}.
\]
The bound \eqref{GZ.bnd} from Guillop\'e-Zworski \cite{GZ:1995a,GZ:1997} gives $N_X(t) \asymp t^2$ for 
more general surfaces with hyperbolic ends, with arbitrary metric inside a compact set, provided the $0$-volume (a regularization of the infinite area) is nonzero.  The bounds apply in particular to smooth geometrically finite hyperbolic surfaces with only one exception: the parabolic cylinder, $\brak{z\mapsto z+1}\backslash \bbH$, which has only a single resonance.

In Borthwick \cite{Borthwick:2012} a sharp constant was identified for the upper bound:
\[
\limsup_{a\to\infty} \frac{2}{a^2}\int_0^a \frac{N_X(t) - N_X(0)}{t^2}\>dt \le \abs{\chi(X)} + \sum \frac{\ell_j}{4},
\]
where $\chi(X)$ is the Euler characteristic and $\{\ell_j\}$ are the boundary lengths of the hyperbolic funnel ends.  
This upper bound implies that the lower bound could also be written with a constant that depends only on $\chi(X)$ and $\{\ell_j\}$.  No Weyl law which would make these asymptotics more precise is known for infinite-area hyperbolic surfaces.  

The exponent of convergence $\delta$ of a Fuchsian group $\Gamma$ is defined by
\[
\delta := \inf \left\{s\ge 0:\> \sum_{T\in\Gamma} e^{-s d(z,Tw)} < \infty\right\}. 
\]
The limit set $L_\Gamma$ of a Fuchsian group $\Gamma$ is the set of accumulation points of an orbit of $\Gamma$, with respect to the Poincar\'e disk topology.   Patterson \cite{Patterson:1976} and Sullivan \cite{Sullivan:1979} famously proved that $\delta$ is the Hausdorff dimension of the limit set,
\[
\delta = \dim L_\Gamma. 
\]
The group $\Gamma$ is cofinite if and only if $\delta = 1$, so the Patterson-Sullivan theory is most interesting in the infinite-area case.  For the elementary cases, which consist of $\bbH$ and the hyperbolic cylinders obtained when $\Gamma$ is cyclic (assuming $X$ is smooth), we have $\delta=0$.  For a non-elementary surface $\delta>0$.

In his original paper, Patterson \cite{Patterson:1976} proved that if $\delta>\tfrac12$, then $\delta(1-\delta)$ is the lowest eigenvalue of $\Delta_X$.  He extended this result in \cite{Patterson:1988}, proving that for all $0< \delta<1$, $R_X(s)$ has a simple pole at $s = \delta$ and no other poles with $\re s \ge \delta$.   The same proofs give also the corresponding result that $\delta$ is the first zero of the zeta function $Z_X(s)$, defined in \eqref{Z.def}.  

Patterson later expanded on this connection between the resonance set and the Selberg zeta function, proving in \cite{Patterson:1989} a trace formula connecting the logarithmic derivative of $Z_X(s)$ to a regularized trace of the resolvent.  (The analogous formula in finite area follows directly from an application of the Selberg trace formula, but the regularization problem is more serious in the infinite-area case.)
Using this formula as a key ingredient, Patterson-Perry \cite{PP:2001} proved the following result for convex co-compact surfaces, which was extended to surfaces with cusps in Borthwick-Judge-Perry \cite{BJP:2005}:
\begin{theorem}
For $X = \Gamma\backslash \bbH$ geometrically finite, the zero set of the zeta function $Z_X(s)$ is the union of the resonance
$\calR_X$ (taken with multiplicities) and the set of points $s = -k$, $k \in \bbN_0$, with multiplicities $-\chi(X)(2k+1)$.  If $X$ has $n_{\rm c}$ cusps then $Z_X(s)$ has poles of order $n_{\rm c}$ at $s \in \tfrac12 - \bbN_0$.
\end{theorem}
\noindent
This connection forms the basis for all of the resonance calculations in this paper; we will describe the methods used to calculating approximate zeros of the Selberg zeta function in \S\ref{zeta.sec}--\ref{comp.sec}.

\section{Zeta function for Schottky groups}\label{zeta.sec}

An argument of Button \cite{Button:1998} shows that any convex co-compact Fuchsian group is a classical Schottky group.
This means that for any smooth, geometrically finite hyperbolic surface $X$ without cusps, there is a classical Schottky group $\Gamma$ such that 
$X \cong \Gamma\backslash \bbH$.   

Suppose the Fuchsian group $\Gamma \in PSL(2,\bbR)$ is a classical Schottky group with $r$ generators.   Such a group is defined by a collection of Euclidean disks $D_1,\dots, D_{2r}$ in $\bbC$ with centers on the real axis and mutually disjoint closures.  If $S_j$ denotes the element of $PSL(2,\bbR)$ that maps the exterior of $D_j$ onto the interior of $D_{j+r}$, then 
$\Gamma$ is freely generated by the transformations $S_1,\dots, S_r$.  
We adopt the standard convention that $S_j$ is defined cyclically for $j\in\bbN$, such that
\[
S_{j+2r} = S_j \quad\text{and}\quad S_{j+r} = S_j^{-1}.
\]
Removing the set of disks $\{D_j\}$ from $\bbH$ yields a fundemantal domain for the action of $\Gamma$, and by triangulating this domain we can easily see that the quotient $X = \Gamma\backslash \bbH$ has Euler characteristic $\chi(X) = 1-r$.

Associated to $\Gamma$ is a Bowen-Series map $B:U\to \bbC\cup\{\infty\}$, where $U := \cup_{j=1}^{2r}{D_j}$, defined simply by
\[
B|_{D_j} = S_j|_{D_j}.
\]
The associated Ruelle transfer operator, acting on the Hilbert space $\calH(U)$ of square-integrable analytic functions, is defined by
\[
(L(s)u)(z) := \sum_{w\in U:\>Bw =z} B'(w)^{-s} u(w).
\]
This operator is trace class on $\calH(U)$ for any $s\in\bbC$, and by Pollicott \cite{Pollicott:1991} we can use it to write the Selberg zeta function as a determinant: 
\begin{theorem}
For all $s\in \bbC$,
\begin{equation}\label{Z.det}
Z_X(s) = \det (I - L(s)).
\end{equation}
\end{theorem}

See \cite[Thm.~5.8]{Borthwick} for details of the proof in the setting of hyperbolic surfaces.  Note that this formula gives a rather direct proof that $Z_X(s)$ extends to an entire function in the convex-co-compact case.  (Conversely, the existence of poles in the zeta function in the case of a surface with cusps shows that we could not extend these methods to that case.)

Our numerical method rests on a algorithm adapted from Jenkinson-Pollicott \cite{JP:2002} and Guillop\'e-Lin-Zworski \cite{GLZ:2004} for computing the coefficients in the power series expansion,
\[
\det (I - z L(s)) = 1 + \sum_{n=1}^\infty d_n(s) z^n.
\]
Using a lemma of Grothendieck \cite[Lemma~2]{JP:2002}, Pollicott-Rocha \cite[Prop~3.6]{PR:1997} prove that for any compact $K \subset \bbC$ there exist constants $C>0$, $c>0$ such that
\begin{equation}\label{dns.est}
\abs{d_n(s)} \le Ce^{-cn^2},\quad\text{for all }s\in K.
\end{equation}
Hence we have a series expansion for the zeta function,
\begin{equation}\label{zeta.sum}
Z_X(s) = 1 + \sum_{n=1}^\infty d_n(s),
\end{equation}
converging uniformly on compact sets in $\bbC$.

The formula for $d_n(s)$ comes from expanding the exponential in 
\[
\det (I - z L(s)) = \exp\left(- \sum_{n=1}^\infty \frac{z^n}{n} \tr L(s)^n \right),
\]
and then collecting powers of $z^n$.  As noted in \cite[\S7.2]{GLZ:2004}, it is convenient to compute $d_n(s)$ as a sum
\begin{equation}\label{dns.def}
d_n(s) := \sum_{k=1}^n b_{n,k}(s),
\end{equation}
where the $b_{n,k}(s)$ are defined recursively by first setting
\[
a_n(s) := - \frac{1}{n} \tr L(s)^n,
\]
and then taking $b_{n,1}(s) := a_n(s)$ and
\[
b_{n,k}(s) := \frac{1}{k} \sum_{m=1}^{n-k-1} b_{n-m,k-1}(s) a_{m}(s),
\]
for $k>1$.

This recursive construction reduces the problem to a computation of $a_n(s)$.  
Under the decomposition $\calH = \oplus_{j=1}^{2r} \calH(D_j)$, the component
\[
L_{ji}: \calH(D_i) \to \calH(D_j), \quad i \ne j+r,
\]
is given by 
\[
L_{ji}u(z) := \left[ (S_i^{-1})'(z)\right]^s u(S_i^{-1}z)\Bigr|_{z\in D_j}
\]
The components of $L(s)^n$ are indexed by
\[
\calW_n := \Bigl\{\sigma \in (\bbZ/2r\bbZ)^n:\> \sigma_{j+1} \ne \sigma_{j} + r\text{ for }j=1,\dots,n-1, \text{and }\sigma_1 \ne \sigma_n+ r\Bigr\},
\]
so that
\[
\tr L(s)^n = \sum_{\sigma\in \calW_n} \tr \left[ L_{\sigma_1\sigma_2}(s) \dots L_{\sigma_n\sigma_1}(s)\right].
\]
If $T_\sigma := S_{\sigma_1}\cdots S_{\sigma_n}$, then the holomorphic fixed point formula (see, e.g.~\cite[Lemma~15.7]{Borthwick}) gives
\[
\tr \left[ L_{\sigma_1\sigma_2}(s) \dots L_{\sigma_n\sigma_1}(s)\right] = \frac{(T_\sigma^{-1})'(w_-)^s}{1 - (T_\sigma^{-1})'(w_-)},
\]
where $w_- \in \bbR\cup\infty$ is the repelling fixed point of $T_\sigma$.
An easy computation shows that $(T^{-1})'(w_-) = e^{-\ell(T_\sigma)}$, so the formula for $a_n(s)$ becomes
\begin{equation}\label{an.sum}
a_n(s) := - \frac{1}{n} \sum_{\sigma\in \calW_n} \frac{e^{-s\ell(T_\sigma)}}{1 - e^{-\ell(T_\sigma)}}.
\end{equation}

From \eqref{an.sum} we can see that the information required to compute $Z_X(s)$ from \eqref{zeta.sum} consists of the set of lengths
\begin{equation}\label{ellWn}
\Lambda_n := \{ \ell(T_\sigma):\>\sigma \in \calW_n \},
\end{equation}
for each $n\in \bbN$.  Note that the multiplicities in this set do not agree with those of the length spectrum itself, because the cyclic permutations of a given index set are listed separately in $\calW_n$.  To explain the precise connection, each element $\ell$ of the length spectrum corresponds to a non-trivial conjugacy class within the group.  For each such conjugacy class there is a minimal word length $n$ needed to write a representative of the class
in terms of the basis $\{S_j\}$.  The length $\ell$ associated to the conjugacy class appears in the set $\Lambda_n$ only for $n$ given by this minimum word length, but with a multiplicity given by the number of representatives of length $n$ within the conjugacy class.  This is the same as the number of distinct cyclic permutations of the word.

\section{Computation}\label{comp.sec}

Jenkinson-Pollicott \cite{JP:2002} proposed using the truncated series,
\begin{equation}\label{ZXN.def}
Z_{X,N}(s) := 1 + \sum_{n=1}^N d_n(s),
\end{equation}
as a numerical method for computing $\delta$, as the convergence turns out to be very rapid near $s = \delta$.  
Guillop\'e-Lin-Zworski \cite{GLZ:2004} adapted this technique to study the growth of the zeta function and the resonance count in strips, for Schottky reflection groups.  In this paper we apply this same approach to the approximation of $Z_{X}(s)$ for convex co-compact Schottky groups, but restricted to groups with two generators.

\subsection{Convergence issues}
In \eqref{an.sum}, we saw that the information required for the calculation of \eqref{ZXN.def} is the collection of sets $\Lambda_n$ for word lengths $1\le n \le N$.   The number of words of each length is
\[
\#\calW_n = (2r-1)^{n} + O(1),
\]
so to calculate $Z_{X,N}(s)$ requires a computation of roughly $(2r-1)^{N+1}$ elements of the length spectrum.  This makes it somewhat impractical to apply the algorithm unless $r=2$.

\begin{figure}
\begin{tabular}{cc}
\begin{overpic}[scale=.7]{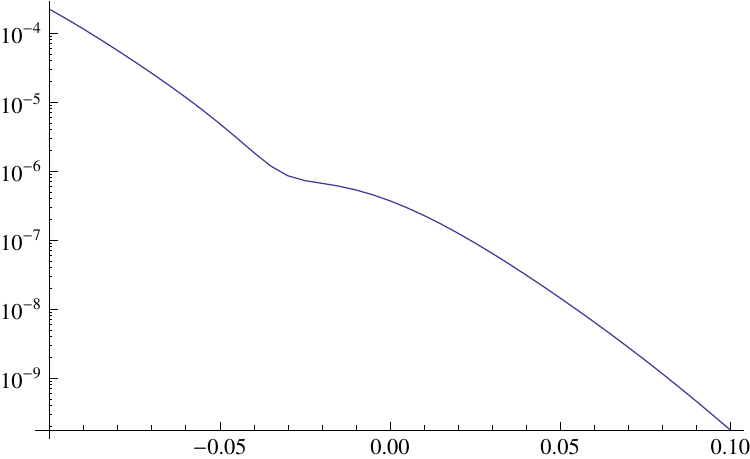}
\put(37,47){$R_{10}(t + 500i)$}
\put(98,6){$t$}
\end{overpic} & 
\begin{overpic}[scale=.7]{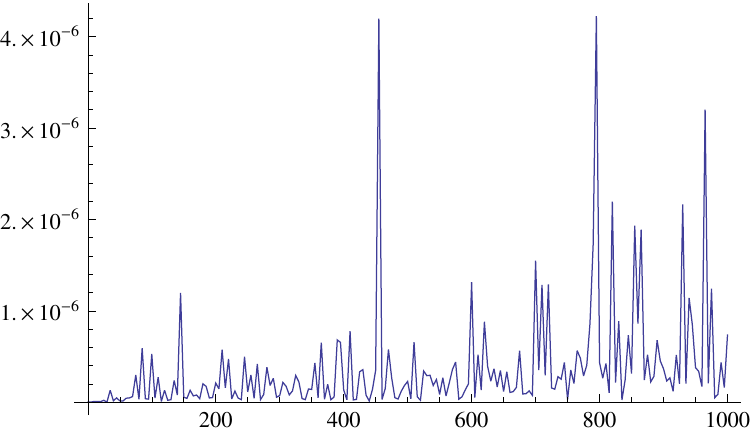}
\put(82,47){$R_{10}(.01+it)$}
\put(99,6){$t$}
\end{overpic}
\end{tabular}
\caption{Behavior of the relative error term $R_n(s) := \abs{d_N(s)/Z_{X,N}(s)}$, for $X(12,14,15)$.}\label{Errors.fig}
\end{figure}

The dependence on $\re s$ presents another serious limitation.  The terms in the sum \eqref{an.sum} for $a_n(s)$ grow exponentially for $\re s < 0$.  
Because of this, the relative errors $d_n(s)/Z_{X,n}(s)$ also grow exponentially as $\re s$ decreases, as illustrated on the right side in 
Figure~\ref{Errors.fig}.  While the computations are still feasible for small negative values of $\re s$ in many cases, the range is quite limited in this direction.  On the other hand, for $\re s \ge 0$ the errors grow quite slowly as $\im s \to \infty$, so we can handle relatively large values of $\im s$.

Another significant factor in the size of error terms is the presence or absence of short geodesics.  If $X$ has many short geodesics then the $a_n(s)$ terms become much larger, and the convergence of \eqref{zeta.sum} can be impractically slow in this case as well.   Longer geodesics correspond in general to smaller values of $\delta$, so this limitation restricts our attention to $\delta$ on the order of $1/10$.

One area where we can gain efficiency is in grouping the elements of $\Lambda_n$ according to multiplicities.
The length spectrum $\calL_X$ of a hyperbolic surface was shown to have unbounded multiplicities by Randol \cite{Randol:1980}.  More recently, Ginzburg-Rudnick \cite{GR:1998} investigated the length spectrum multiplicity problem for the case of a free group $\Gamma$ on two generators $A, B$, e.g.~a classical Schottky group with $r=2$.   As mentioned above, the multiplicity problem in $\Lambda_n$ is not the true length multiplicity problem, but is related in a direct way.

Define the \emph{length class} of an element $\sigma \in \calW_n$ by
\[
[\sigma] := \bigl\{\sigma'\in\calW_n:\> \ell(T_{\sigma'}) = \ell(T_{\sigma})\text{ for any choice of generators}\bigr\}.
\]
The full multiplicity of $\ell(T_\sigma)$ in $\Lambda_n$ will be the size of the class $\#[\sigma]$.  
The operations on $\sigma$ that clearly preserve the length class are:
\begin{enumerate}
\item  Cyclic permutation:  $(\sigma_1,\dots,\sigma_n)\mapsto (\sigma_2,\dots,\sigma_n,\sigma_1)$ 
\item  Inverse: $T_\sigma \mapsto T_\sigma^{-1}$.
\item  Reverse:  $(\sigma_1,\dots,\sigma_n) \mapsto (\sigma_n,\dots,\sigma_1)$.
\item  Decompose $T_\sigma$ as a word $w(U,V)$ and then reverse the order by
\[w(U,V) \mapsto w(U^{-1}, V^{-1})^{-1}.\]
\end{enumerate}
A typical word of length $n$ generates a length class of multiplicity $4n$, consisting $n$ permutations each for the word, its reverse, its inverse, and the reverse of its inverse.  Words with greater  
multiplicity can be formed by applying these symmetry operations to composite words.

According to \cite[Conjecture~1]{GR:1998}, this should be a complete list of operations that preserve the length class.  
This is not proven, and in any case the operations do not always lead to distinct elements of $\calW_n$, so it remains a difficult combinatorial problem to understand the structure of the length classes for a general value of $n$.

As part of the strategy for efficient computation of the functions $a_n(s)$, we first compute a multiplicity table breaking down $\calW_n$ into length classes.  This is done simply by computing lengths for all elements of $\calW_n$ for a few different sets of generators and then correlating the results.
Since the groups are freely generated, the structure of the length classes depends only on the number of generators.  Table~\ref{Lnsize.tab} shows a comparison of the number of length classes to the number of elements of $\Lambda_n$.

\begin{table}
\begin{center}
\setlength{\extrarowheight}{3pt}
\begin{tabular}{|c|c|c|}
\hline $n$  & $\#\Lambda_n$ & $\#\{\text{length classes}\subset \Lambda_n\}$  \\
\hline 1 & 4 & 2 \\
\hline 2 & 12 & 4 \\
\hline 3 & 28 & 6 \\
\hline 4 & 84 & 13 \\
\hline 5 & 244 & 22 \\
\hline 6 & 732 & 52 \\
\hline 7 & 2188 & 106 \\
\hline 8 & 6564 & 266 \\
\hline 9 & 19684 & 626 \\
\hline 10 & 59052 & 1632 \\
\hline 11 & 177148 & 4218 \\
\hline 12 & 531444 & 11471 \\
\hline \end{tabular}
\end{center}
\caption{Size of $\Lambda_n$ versus number of distinct lengths.}\label{Lnsize.tab}
\end{table}

\subsection{Locating zeros}

\begin{figure}
\begin{center}
\begin{overpic}[scale=.6]{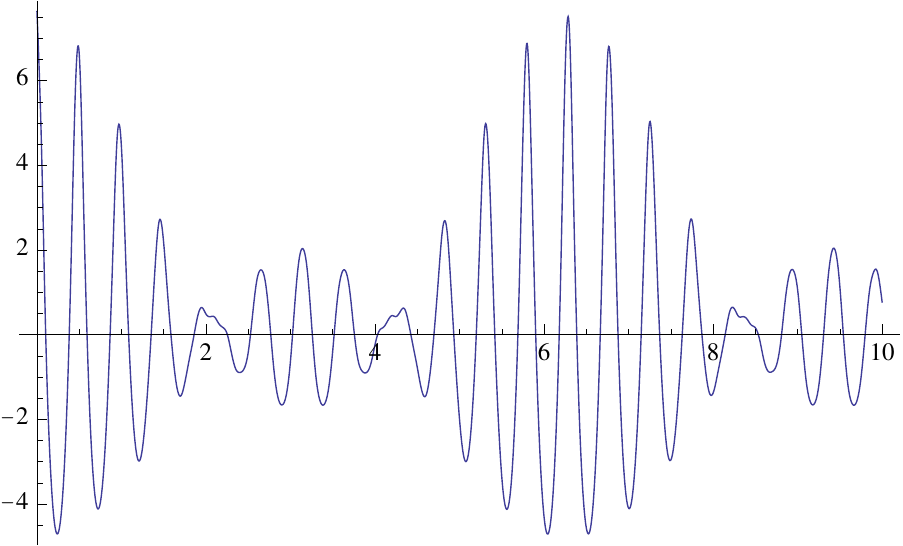}
\end{overpic}
\end{center}
\caption{Sample oscillations in $\re Z_X'/Z_X(s)$, with $s = 2\delta + it$, for the surface $X(12,13,14)$ for which $\delta \approx 0.1068$.}\label{logZp.fig}
\end{figure}

The argument principle is the obvious technique for locating or counting zeros of an analytic function.  The derivative $Z_X'(s)$ can be computed as in \eqref{zeta.sum},
simply be replacing $a_n(s)$ with the derivative,
\[
a_n'(s) =  \frac{1}{n} \sum_{\sigma\in \calW_n} \frac{\ell(T_\sigma) e^{-s\ell(T_\sigma)}}{1 - e^{-\ell(T_\sigma)}},
\]
and then running the same recursive algorithm used to compute $d_n(s)$.
However the logarithmic derivative $Z_{X}'/Z_{X}(s)$ is highly oscillatory, even away from the zeros, as illustrated in Figure~\ref{logZp.fig}.  For $\re s$ closer to $0$ these oscillations become much more severe, so that numerical integration of $Z_X'/Z_X(s)$ over long intervals is extremely slow and prone to errors.   
Another issue with this approach is that the numerical routines often fail when poles are close to the path.  In regions with a high density of poles, accurate estimation of the integrals requires adjustments to the paths to avoid poles.

We have found that a much faster and more stable approach to the application of the argument principle is 
subdivide the sample region into a rectangular grid.  Along each short horizontal or vertical segment in the grid, 
we sample $Z_X(s)$ at a finely-spaced sequence of points $\{s_j\}$,
and then compute 
\[
\varDelta \arg Z_X = \sum_{j=1}^m \arg \frac{Z_X(s_{j+1})}{Z_X(s_j)}.
\]
We can efficiently compute the resonance counts in a rectangular bin array by first using this sampling approach to obtain $\varDelta \arg Z_X$ on all of the horizontal and vertical edge segments.   For the sampling we specify a minimum spacing, typically 0.01, as well as a minimum number of samples per edge, but otherwise the samples are spaced uniformly.  The stability of the method is easily checked by counting zeros in the same region with varying bin shapes and sizes.  
The method yields correct results even when zeros lie extremely close to the edges of the bins,
so we are able to avoid the issue of adjusting paths to avoid poles.

For the resonance plots shown in Figure~\ref{ResPlot101.fig} and below, we start by producing a bin count with a very fine mesh, so that each bin is roughly the size of a pixel in the final plot.  Then for each bin with a positive count, we locate the corresponding zero using a standard root-finding algorithm, with seed point furnished by the center of the bin.  Although bins with multiple resonances are possibly underrepresented by this technique, this has no impact on the plot because the resulting points would be superimposed and indistinguishable.

\section{Resonance plots}\label{resplot.sec}

Since we restrict our attention to Schottky groups with $r=2$ generators, we can only produce hyperbolic surfaces of Euler characteristic $\chi = -1$.  There are two possible topologies: a surface of genus zero with three funnels or a surface of genus one with one funnel.

\begin{figure}
\begin{center}
\begin{overpic}[scale=.9]{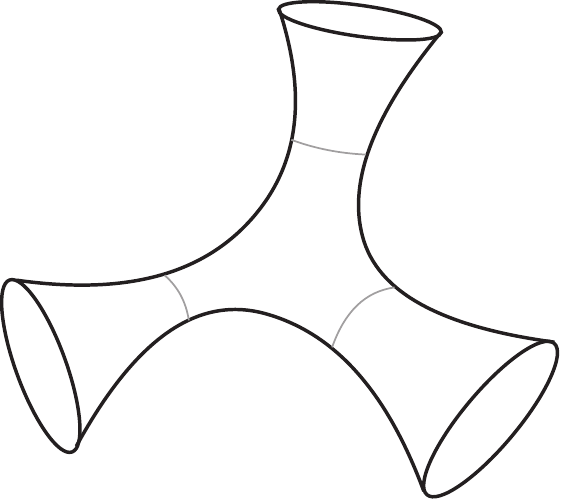}
\put(26,43){$\ell_1$}
\put(45,65){$\ell_2$}
\put(72,38){$\ell_3$}
\end{overpic}
\end{center}
\caption{Schematic view of a 3-funnel surface $X(\ell_1,\ell_2,\ell_3)$.}\label{trif.fig}
\end{figure}

\subsection{3-funnel surfaces}
A hyperbolic pair of pants is bounded by 3 simple closed geodesics.  There is a unique way to attach funnels to these boundary circles, producing a 3-funnel surface topologically equivalent to a thrice punctured sphere.  We denote such a surface by $X(\ell_1,\ell_2,\ell_3)$, where the $\ell_j$ are boundary lengths of the original pair of pants.  These three real parameters characterize the surface completely up to isometry.  

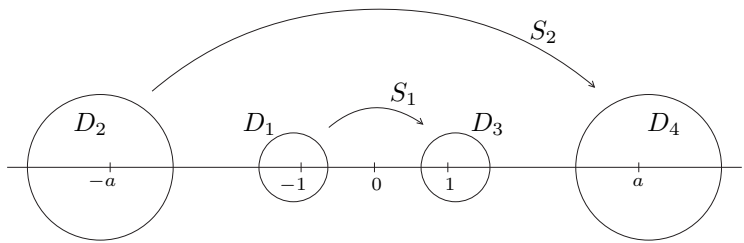
\begin{figure}
\begin{center}
\begin{overpic}[scale=.6]{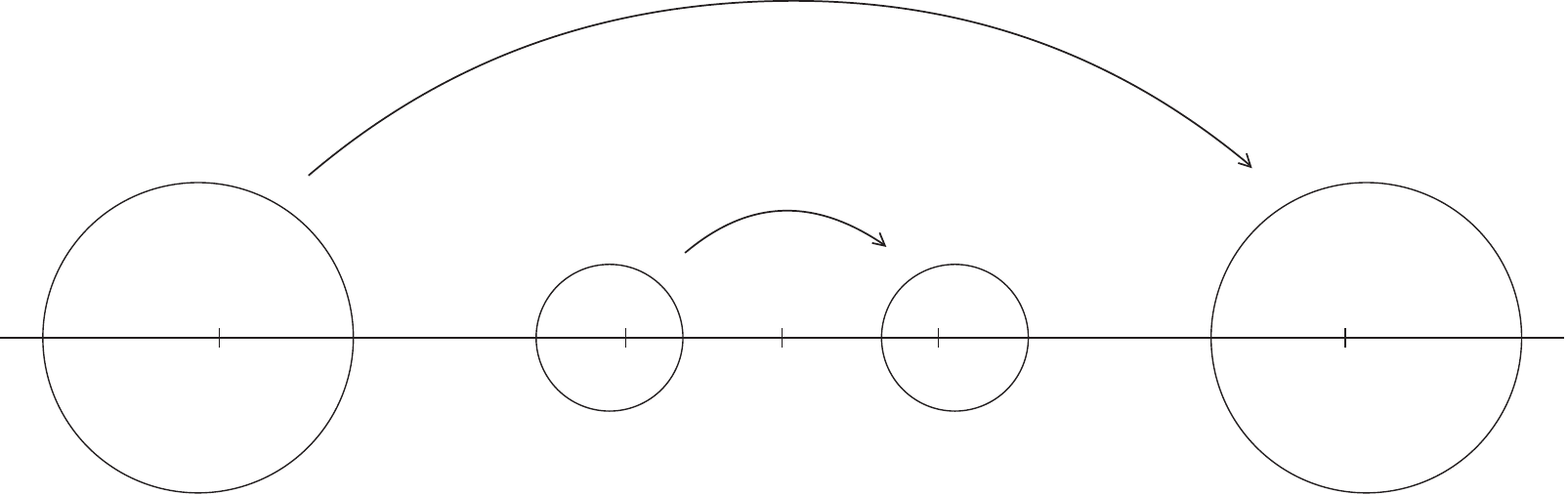}
\put(9,15){$D_2$}
\put(32,15){$D_1$}
\put(63,15){$D_3$}
\put(87,15){$D_4$}
\put(11,7.5){$\scriptstyle -a$}
\put(37,7){$\scriptstyle -1$}
\put(49.5,7){$\scriptstyle 0$}
\put(59.5,7){$\scriptstyle 1$}
\put(85,7.5){$\scriptstyle a$}
\put(52,19){$S_1$}
\put(71,27.5){$S_2$}
\end{overpic} 
\end{center}
\caption{Schottky circles and generators for a 3-funnel surface.}\label{circles1.fig}
\end{figure}

To associate a Schottky group to $\ell_1,\ell_2,\ell_3$, we set
\[
S_1 := \begin{pmatrix}\cosh(\ell_1/2) & \sinh(\ell_1/2) \\  \sinh(\ell_1/2)& \cosh(\ell_1/2) \end{pmatrix},\qquad
S_2 := \begin{pmatrix}\cosh(\ell_2/2) & a \sinh(\ell_2/2) \\  a^{-1} \sinh(\ell_2/2)& \cosh(\ell_2/2) \end{pmatrix},
\]
with the parameter $a>0$ chosen by solving $\tr(S_1 S_2^{-1}) = 2 (\cosh \ell_3/2)$.  In this scheme, the circles $D_1, D_3$ are centered at $\pm 1$, and $D_2, D_4$ at $\pm a$, as shown in Figure~\ref{circles1.fig}.

\begin{figure}
\begin{tabular}{cc}
\begin{overpic}[scale=.7]{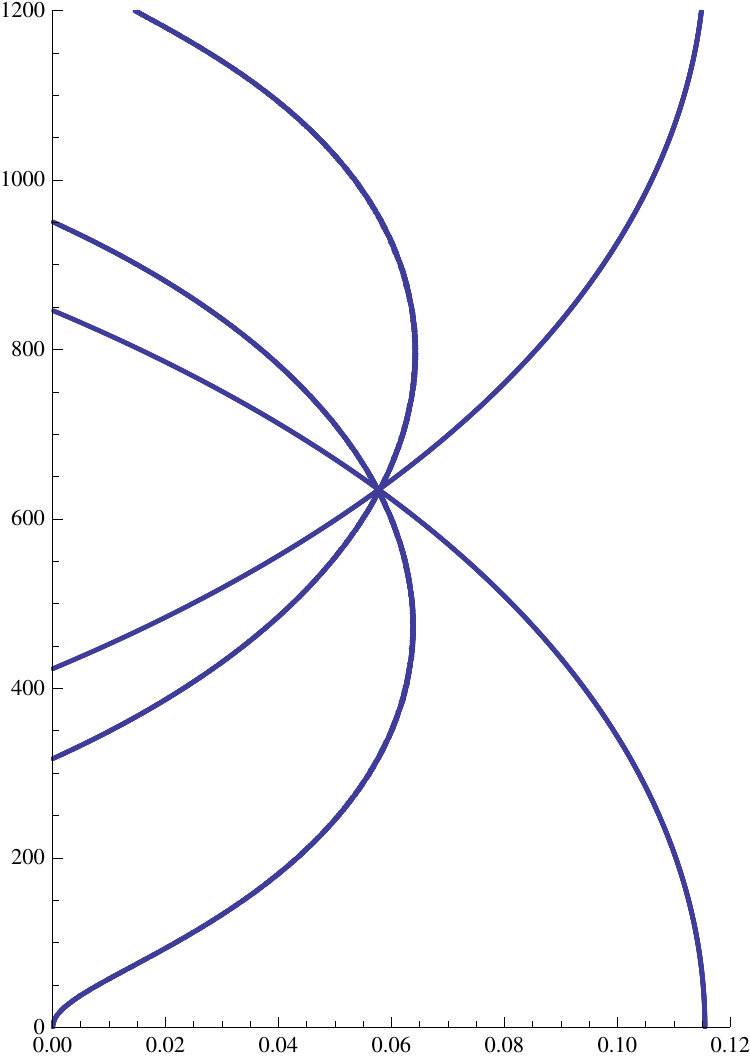}
\put(25,-3){$\scriptstyle X(12,12,12)$}
\end{overpic}  &
\begin{overpic}[scale=.7]{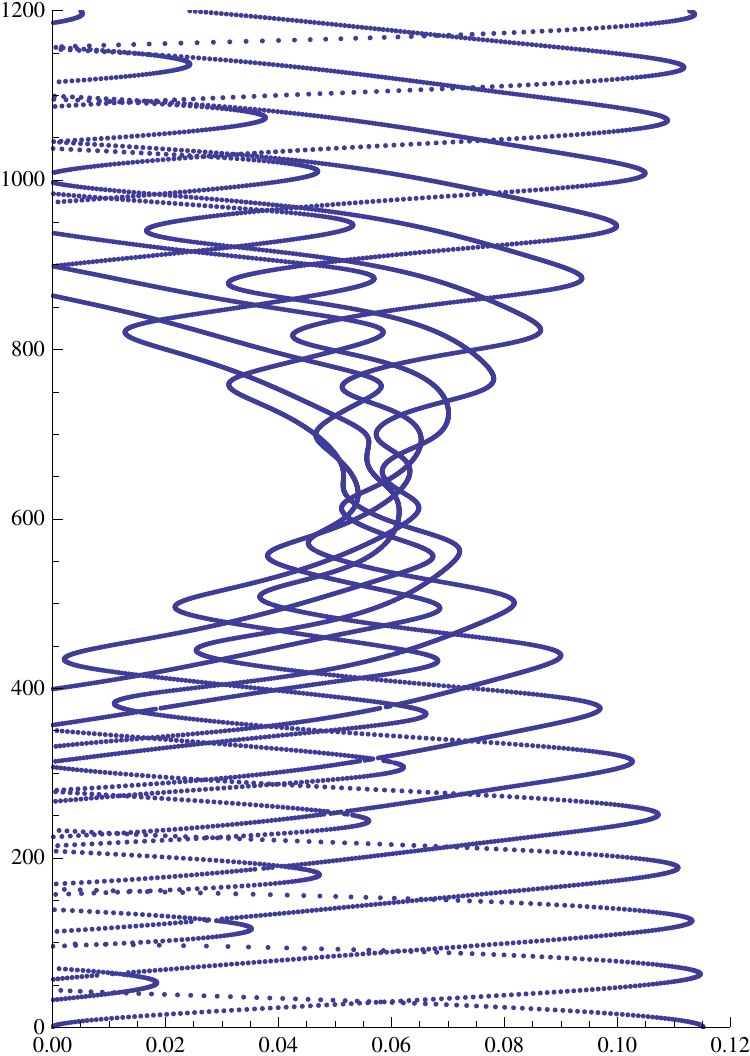}
\put(25,-3){$\scriptstyle X(12,12,12.1)$}
\end{overpic} 
\end{tabular}
\caption{Resonances of 3-funnel surfaces.}\label{ResPlot121212.fig}
\end{figure}

Some sample resonance plots for 3-funnel surfaces were shown in Figure~\ref{ResPlot101.fig}.
In the symmetric case, $X(\ell,\ell,\ell)$, the resonances are concentrated very strongly a few distinctive resonance curves.  
Along some of these, for example the curve that starts at the origin, the multiplicity of each resonance is 2.
The plot on the left in Figure~\ref{ResPlot121212.fig} shows another example of this structure.  These resonance curves are strongly reminiscent of the resonance families seen in $n$-disk scattering systems, as seen, for example, in the plots of numerically calculated resonances in \cite{GR:1992, WH:1998}.  This similarity is perhaps not surprising, as the classical dynamics of a 3-disk scattering system bear a strong resemblance to the dynamics of the symmetric surface $X(\ell,\ell,\ell)$.  
However, for hyperbolic surfaces the existence of these patterns was not previously predicted.  We believe the origin of the resonance chains can be explained through a symmetry factorization of the zeta function (work in progress).

One very interesting feature, in view of the decay rate conjecture to be discussed in \S\ref{decay.sec}, is the fact that the hubs where resonance curves meet occur on the line $\re s \approx \delta/2$.  This seems to hold true for all of the $X(\ell,\ell,\ell)$ cases.  The hubs do start to lose some coherence as $\im s$ increases, as seen in the plot on the left in Figure~\ref{ResPlot101.fig}.

\begin{figure}
\begin{tabular}{ccc}
\begin{overpic}[scale=.7]{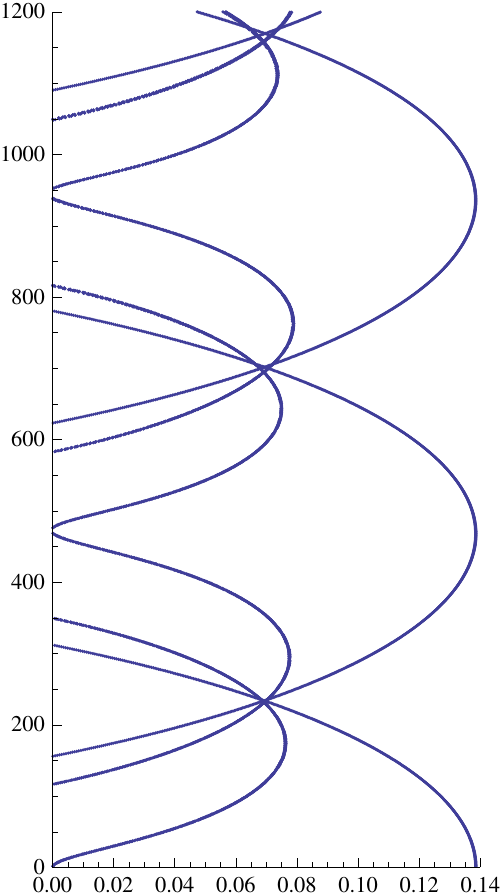}
\put(12,-4){$\scriptstyle X(10,10,10)$}
\end{overpic} &
\begin{overpic}[scale=.7]{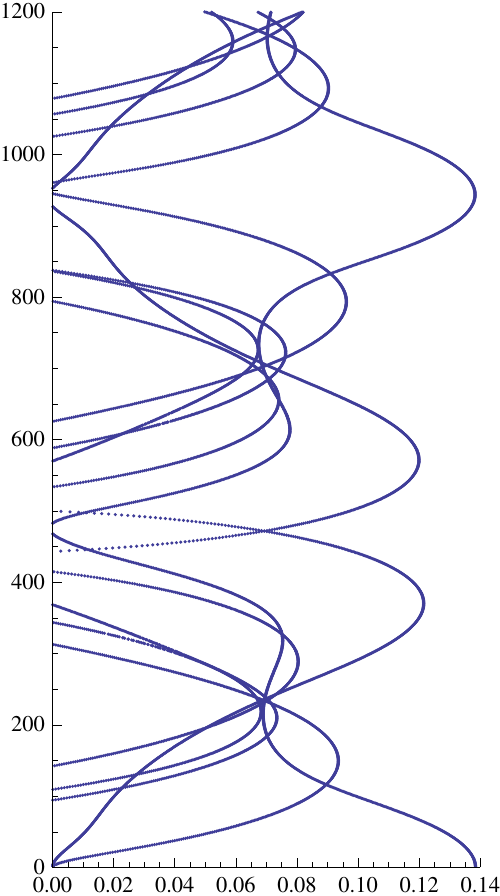}
\put(12,-4){$\scriptstyle X(10,10,10.02)$ }
\end{overpic} &
\begin{overpic}[scale=.7]{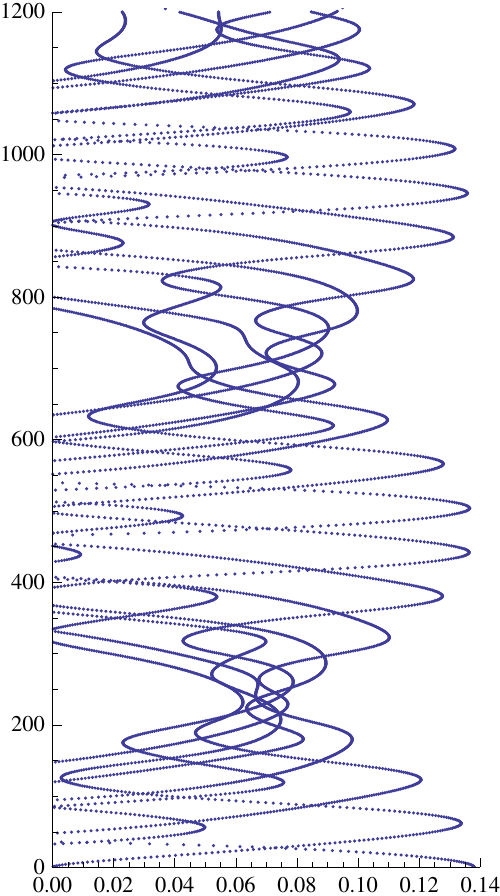}
\put(12,-4){$\scriptstyle X(10,10,10.1)$}
\end{overpic} \\
&& \\
\begin{overpic}[scale=.7]{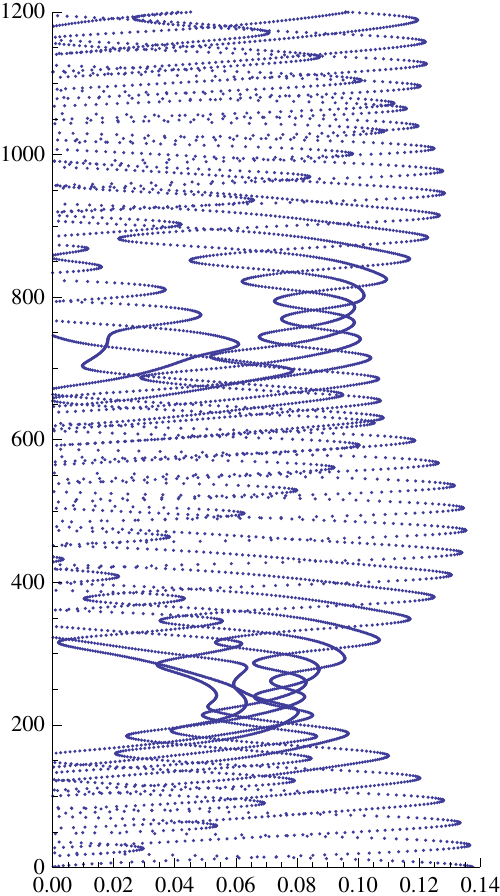}
\put(12,-4){$\scriptstyle X(10,10,10.2)$}
\end{overpic} &
\begin{overpic}[scale=.7]{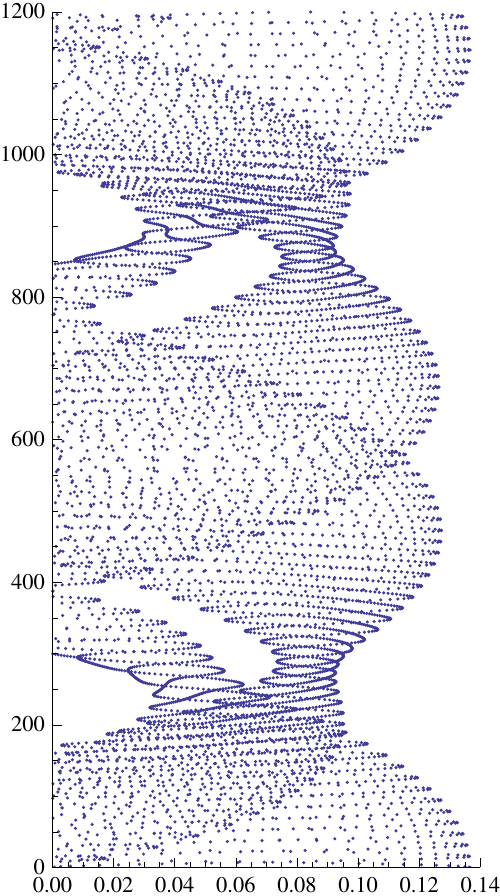}
\put(12,-4){$\scriptstyle X(10,10,10.4)$}
\end{overpic} &
\begin{overpic}[scale=.7]{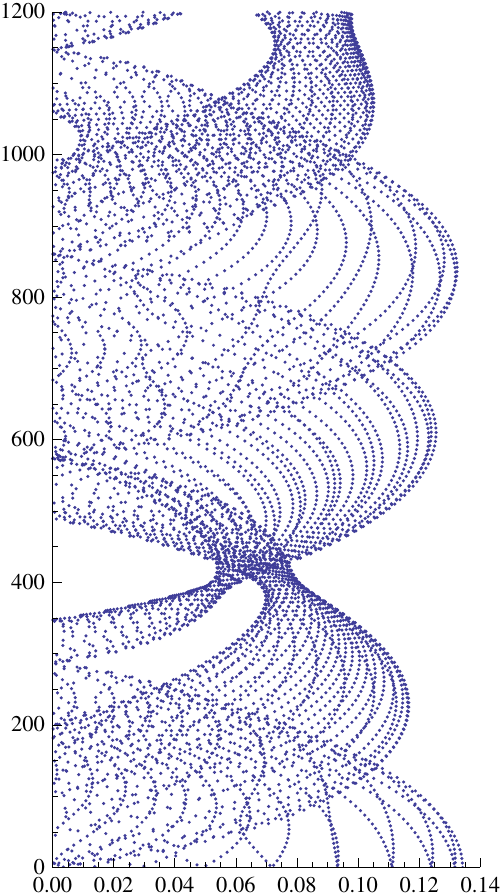}
\put(12,-4){$\scriptstyle X(10,10,11)$}
\end{overpic}
\end{tabular}
\caption{Evolution of resonance patterns away from the symmetric case.}\label{Res10perturbed.fig}
\end{figure}

When we perturb slightly from the symmetric case, as shown on the right in Figures~\ref{ResPlot101.fig} or \ref{ResPlot121212.fig}, we see smaller-scale oscillations of the resonance curves.  The tightly concentrated hubs near $\re s = \delta/2$ remain, but are somewhat broader.
Figure~\ref{Res10perturbed.fig} traces the evolution of the patterns as the case shown in Figure~\ref{ResPlot101.fig}
is further perturbed.

The spacing of the resonances within each resonance curve is extremely consistent.  For example, consider the bottom left arc of the $X(12,12,12)$ plot in Figure~\ref{ResPlot121212.fig}.  Over the first thousand points (from $\im s = 0$ to $\im s \approx 523$) the distance between resonances decreases monotonically from $0.5234905$ to $0.5234901$.  Similarly, for the corresponding arc for $X(10,10,10)$, the spacing decreases from $0.627894$ at an extremely slow rate.  It is interesting to note that $2\pi/12 = 0.523599...$ while $2\pi/10 = .628319...$.  This approximate (and never exact) spacing of $2\pi/\ell$ appears in all of the $X(\ell,\ell,\ell)$ cases, and is stable under small perturbations.

\begin{figure}
\begin{tabular}{cc}
\begin{overpic}[scale=.7]{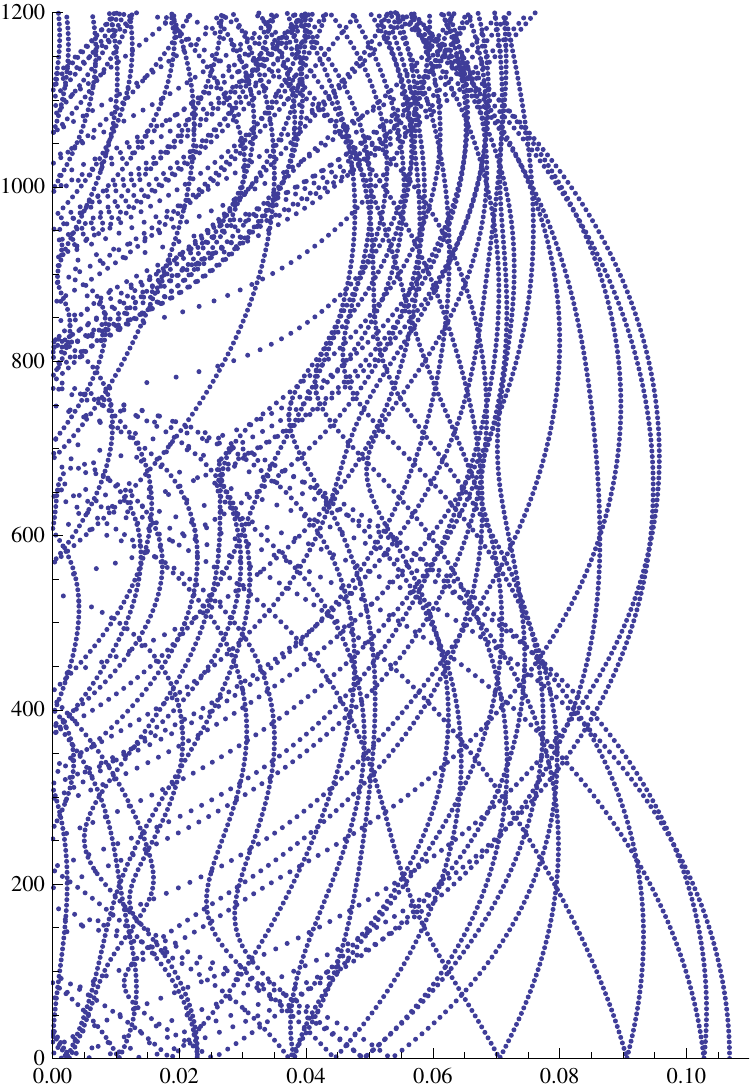}
\put(25,-3){$\scriptstyle X(12,13,14)$}
\end{overpic}  &
\begin{overpic}[scale=.7]{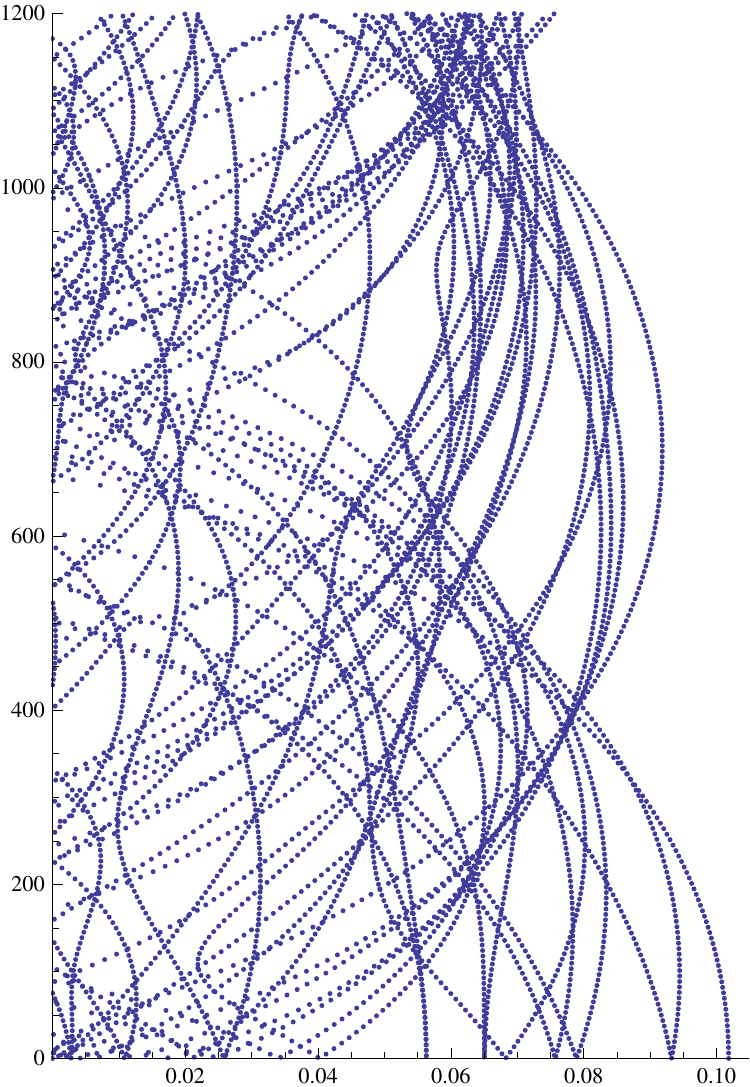}
\put(25,-3){$\scriptstyle X(12,14,15)$}
\end{overpic} 
\end{tabular}
\caption{Resonances of 3-funnel surfaces.}\label{ResPlot121415.fig}
\end{figure}

Further away from the symmetric case, the pattern of resonance curves changes significantly, as shown in Figure~\ref{ResPlot121415.fig}.  
The curves visible on the right in Figure~\ref{ResPlot121212.fig} actually flatten out very quickly as we move away from the symmetric case, and the curves visible in Figure~\ref{ResPlot121415.fig} are actually new patterns of coherence that emerge between points in these flattened curves.  
One can see the onset of this coherence in the final two plots of Figure~\ref{Res10perturbed.fig}.

The spacing along these curves is also very consistent.  For both surfaces shown shown in Figure~\ref{ResPlot121415.fig}, the spacing along each curve is approximately $6.282$.  

\begin{figure}
\begin{center}
\begin{overpic}{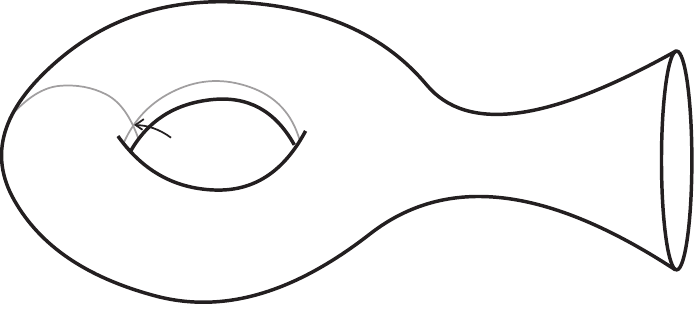}
\put(-1,31){$\ell_1$}
\put(40,30){$\ell_2$}
\put(25,23){$\phi$}
\end{overpic}
\end{center}
\caption{Schematic view of a funneled torus $Y(\ell_1,\ell_2,\phi)$.}\label{torf.fig}
\end{figure}

\subsection{Funneled tori}
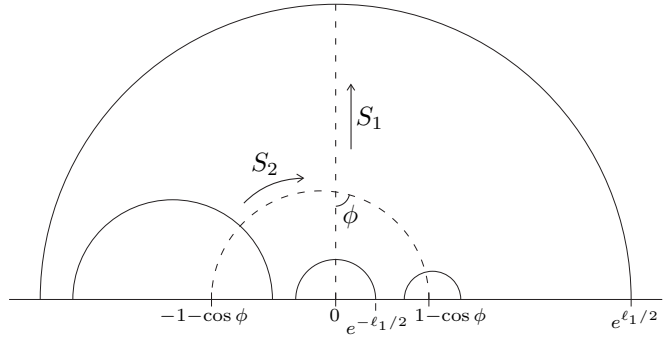
\begin{figure}
\begin{center}
\begin{overpic}[scale=.8]{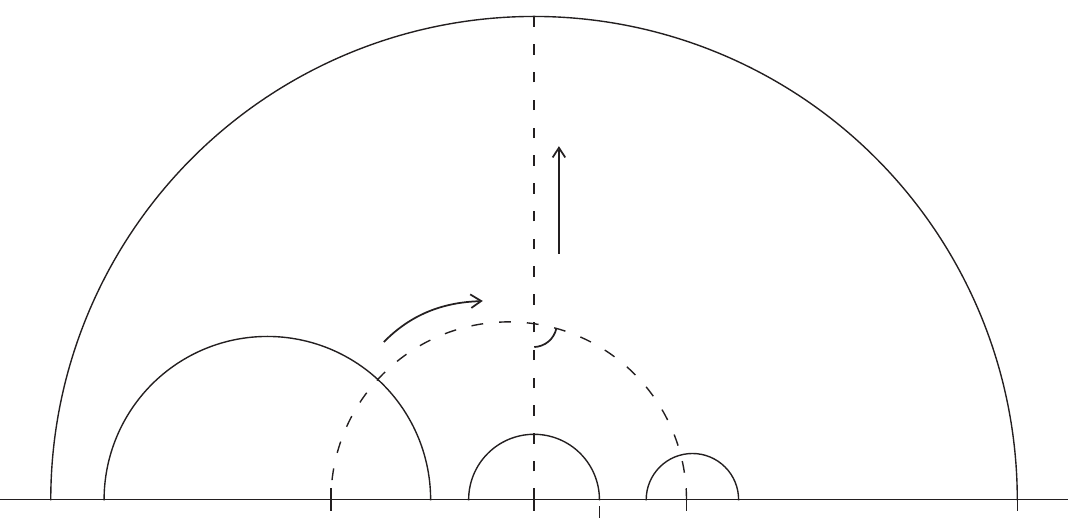}
\put(37,22){$S_2$}
\put(53,29){$S_1$}
\put(51,14){$\phi$}
\put(23,-1){$\scriptstyle -1-\cos\phi$}
\put(49,-1){$\scriptstyle 0$}
\put(51.5,-3){$\scriptstyle e^{-\ell_1/2}$}
\put(62,-1){$\scriptstyle 1-\cos\phi$}
\put(92,-2){$\scriptstyle e^{\ell_1/2}$}
\end{overpic}
\end{center}
\caption{Schottky circles and generators for a funneled torus.}\label{circles2.fig}
\end{figure}

A suitable parametrization of the moduli space for a hyperbolic surface of genus one with one funnel was worked out in Buser-Semmler \cite{BS:1988}.   
The parameters consist of two lengths $\ell_1, \ell_2$, and an  angle $\phi \in (0,\pi)$, and we will denote the corresponding surface as 
$Y(\ell_1,\ell_2,\phi)$.  The corresponding generators in $PSL(2,\bbR)$ are defined by
\[
S_1 = \begin{pmatrix} e^{\ell_1/2} & 0 \\  0 & e^{-\ell_2/2} \end{pmatrix},
\]
\[
S_2 = \begin{pmatrix}\cosh \ell_2/2 - \cos\phi \>\sinh \ell_2/2 & \sin^2\phi \>\sinh \ell_2/2 \\ 
\sinh \ell_2/2 & \cosh \ell_2/2 + \cos\phi \>\sinh \ell_2/2  \end{pmatrix}.
\]
The length parameters are the lengths $\ell_j = \ell(S_j)$ for $j=1,2$ of two simple closed geodesics, with $\phi$ is the angle between them,
as shown in Figure~\ref{torf.fig}.  If $\ell_1,\ell_2$ are sufficiently large and $\phi$ is sufficiently close to $\pi/2$, then $S_1$ and $S_2$ will be Schottky generators, with circles as shown in Figure~\ref{circles2.fig}.

\begin{figure}
\begin{tabular}{cc}
\begin{overpic}[scale=.72]{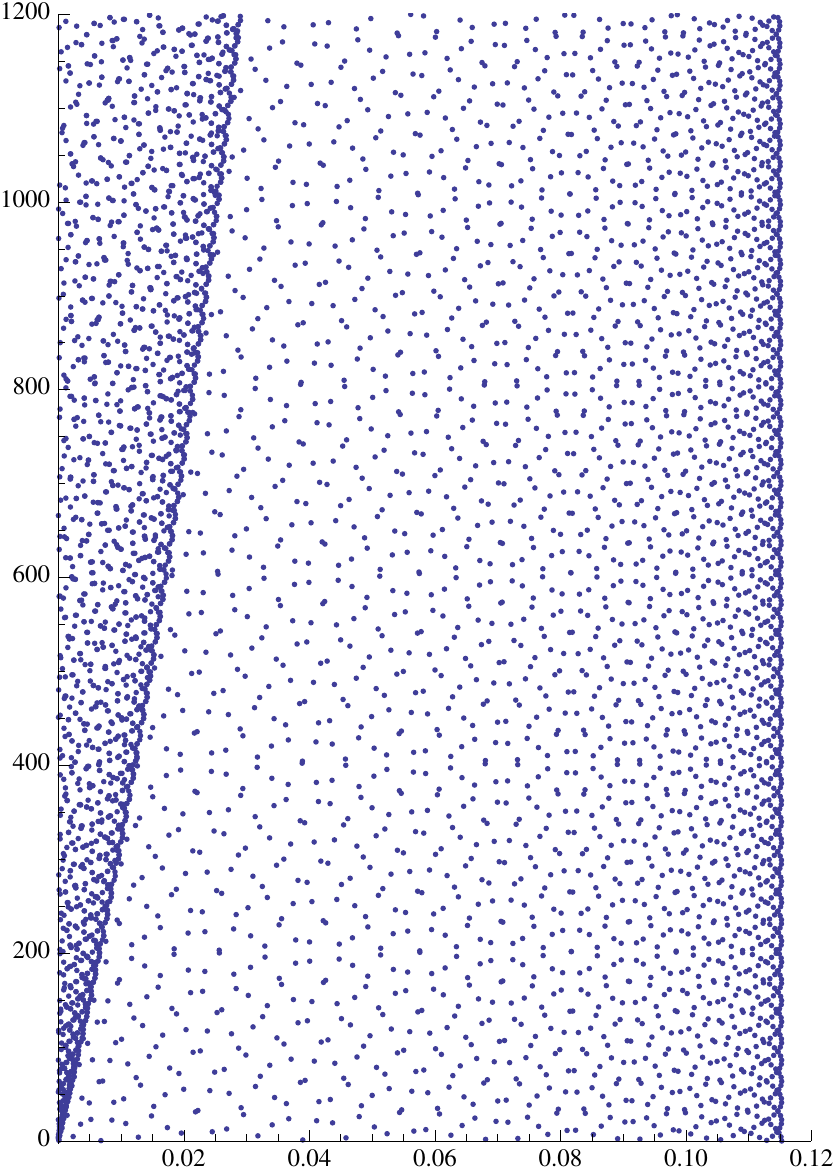}
\put(25,-3){$\scriptstyle Y(10,10,\pi/2)$}
\end{overpic} &
\begin{overpic}[scale=.72]{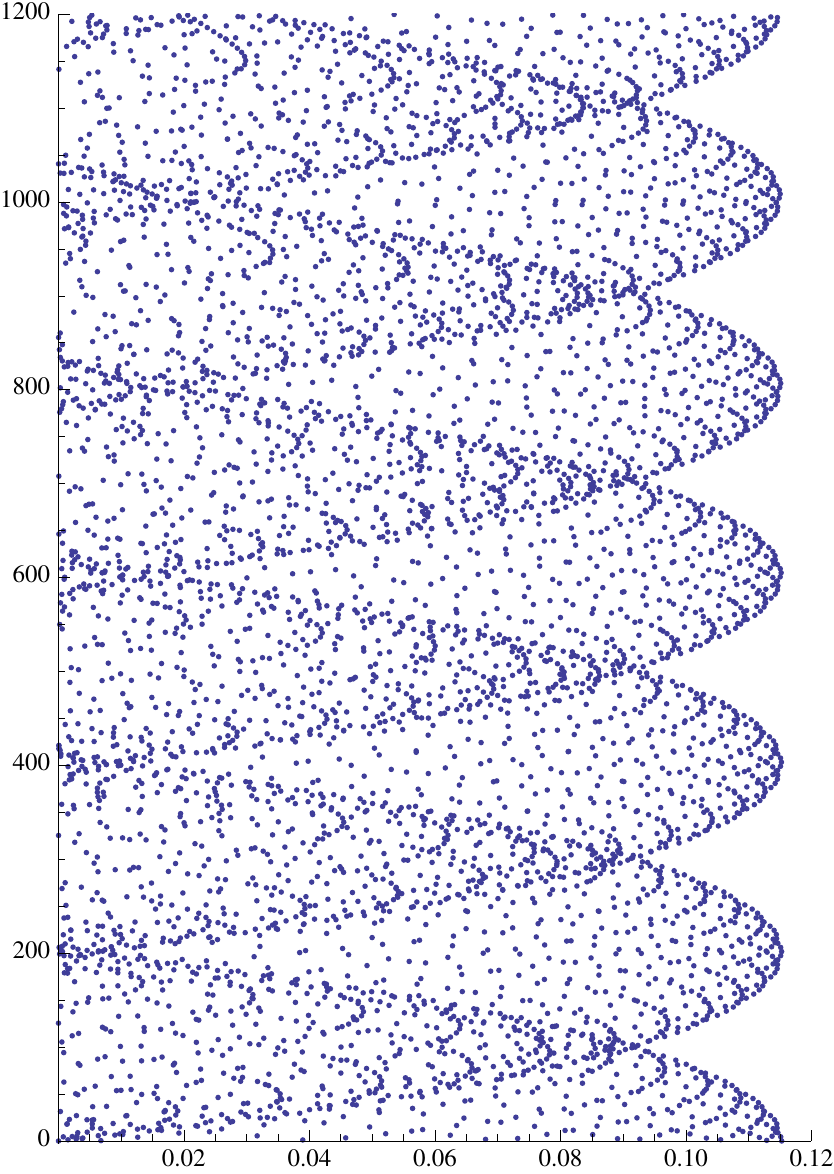}
\put(25,-3){$\scriptstyle Y(10,10,\pi/2.01)$}
\end{overpic}
\end{tabular}
\caption{Resonances of funneled tori.}\label{ResPlotTor.fig}
\end{figure}

\begin{figure}
\begin{center}
\begin{overpic}[scale=.9]{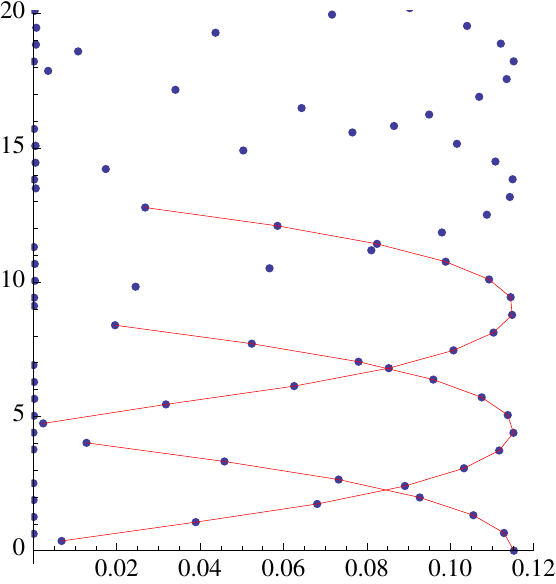}
\end{overpic}
\end{center}
\caption{Resonance curves for $Y(10,10,\pi/2)$.}\label{ResCurve1010.fig}
\end{figure}

A portion of the resonance sets for two funneled tori are shown in Figure~\ref{ResPlotTor.fig}.  The pattern that appears in the plot for $Y(10,10,\pi/2)$ on the left is the result of overlapping resonance curves similar to those seen for the 3-funnel surfaces.  The curves is clearly resolved if we plot on a smaller scale, as illustrated in Figure~\ref{ResCurve1010.fig}.

\begin{figure}
\begin{center}
\begin{overpic}[scale=.6]{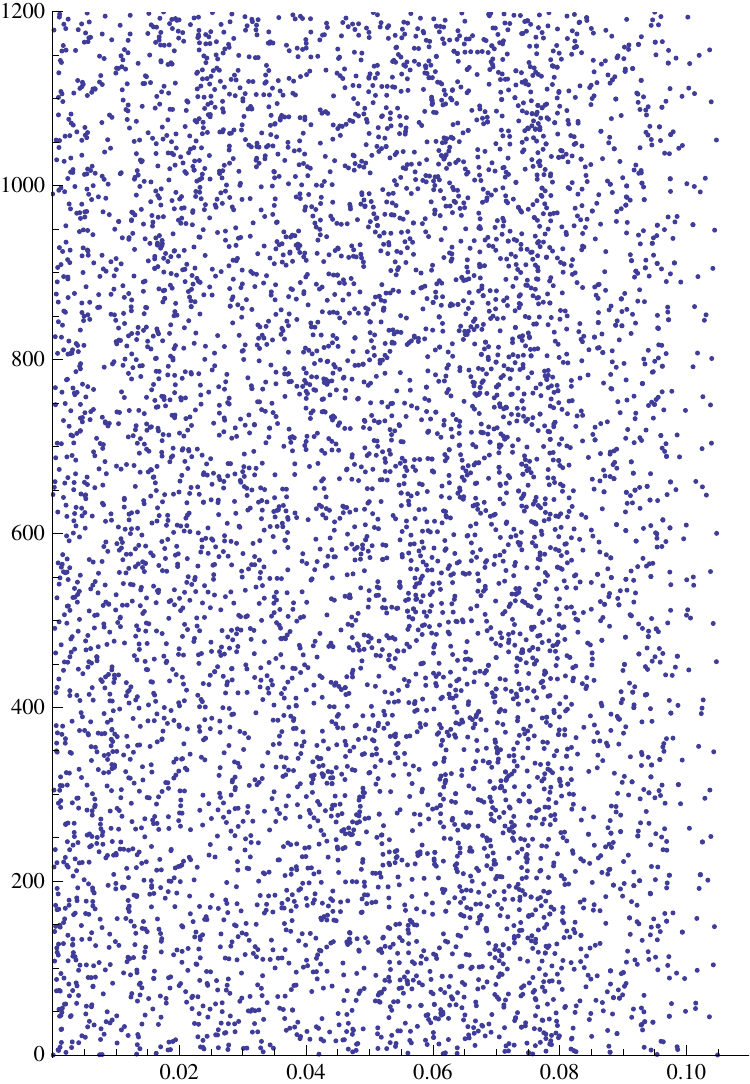}
\end{overpic}
\end{center}
\caption{Resonance plot for $Y(10,12,\pi/2.5)$.}\label{ResPlot1012.fig}
\end{figure}

For funneled tori with parameters further away from the symmetric case, the resonance plots appear incoherent, as shown in Figure~\ref{ResPlot1012.fig}.  The resonance curves do persist if we perturb only very slightly from the symmetric case, as the plot on the right in Figure~\ref{ResPlotTor.fig} illustrates.  However, the pattern changes quite rapidly as the parameters vary.  Presumably patterns of resonance curves are present even in plots such as Figure~\ref{ResPlot1012.fig}, but overlapping to such an extent that they are difficult to distinguish.

\section{Fractal Weyl conjecture}\label{frweyl.sec}

The work of Sj\"ostrand \cite{Sjostrand:1990} on semiclassical bounds for resonance counting has led to a general expectation for chaotic scattering systems that the number of resonances near the continuous spectrum should satisfy a power law with exponent equal to half of the dimension of the classical trapped set.  Recently a large number of theoretical \cite{Zworski:1999, GLZ:2004, NZ:2005, SZ:2007, NR:2007, NZ:2007, Nonn:2008, DD:2012}, numerical \cite{SFPB:2000, Lin:2002, LZ:2002, LSZ:2003, SWM:2009} and experimental \cite{LRPS:1999, LVPRS:2000, PWBKSZ:2012, BWPSKZ:2012, KMBK:2013} studies have appeared in support of this conjectural `fractal Weyl law'.

For a geometrically finite hyperbolic surface, $X = \Gamma\backslash \bbH$, a geodesic is trapped if and only if its lift to $\bbH$ has endpoints in the limit set $L_\Gamma$ introduced in \S\ref{prelim.sec}.  Since $\dim L_\Gamma = \delta$, it follows by a straightfoward calculation that the dimension of the trapped set (as a subset of  $T^*X$) is $2+2\delta$.  

Let us define the resonance counting function for a vertical strip as
\[
N_X(a_0,a_1;t) := \#\Bigl\{\zeta\in\calR_X:\> a_0 \le \re \zeta \le a_1, 0 \le \abs{\im\zeta} \le t\Bigr\}.
\]
For the general case of open chaotic systems the fractal Weyl conjecture was first made formally in Lu-Sridhar-Zworski \cite{LSZ:2003}.  In our context the conjecture could be paraphrased as:    
\begin{conj}[Fractal Weyl law]
For a geometrically finite hyperbolic surface $X$, 
\begin{equation}\label{weyl.conj}
N_X(\sigma,\delta;t) \asymp t^{1+\delta},
\end{equation}
provided that $\sigma$ is sufficiently small.
\end{conj}
\noindent
The conjecture is sometimes stated in a stronger form, as an exact asymptotic, $N_X(\sigma,\delta;t) \sim At^{1+\delta}$ for some constant $A$.  For finite-area surfaces, since $\delta =1$ and all resonances are contained in a vertical strip, the fractal Weyl law is simply the Weyl law \eqref{finite.weyl}.

Guillop\'e-Lin-Zworski \cite{GLZ:2004}, improving on an earlier result of Zworski \cite{Zworski:1999}, proved a stronger, localized version of the conjectured upper bound: for any $\sigma$,
\begin{equation}\label{glz.bound}
N_X(\sigma,\delta;t+1) - N_X(\sigma,\delta;t) = O_\sigma(t^{\delta}).
\end{equation}
Note that this is consistent with a remainder term of $O(t^{\delta})$ for the upper bound in \eqref{weyl.conj}.   This estimate \eqref{glz.bound} apples to hyperbolic Schottky manifolds in any dimension,  and has been extended to general asymptotically hyperbolic manifolds by Datchev-Dyatlov \cite{DD:2012}, with $\delta$ reinterpreted in terms of the upper Minkowski dimension of the trapped set.  For conformally compact hyperbolic surfaces, Dyatlov-Guillarmou \cite{DG:2012} have established a fractal bound on the remainder term in the Weyl asymptotic for the scattering phase, also consistent with the fractal Weyl law.

The best known lower bound comes from the trace formula and is considerably weaker than the conjectured growth rate.  Guillop\'e-Zworski \cite{GZ:1999} proved that for any 
$c > \tfrac32$,
\[
N_X(-c,\delta;t) \ne O\left(t^{1-\frac{2}{2c+1}}\right).
\]
Jakobson-Naud \cite{JN:2010} proved a lower bound that shows explicit dependence on $\delta$, but the growth rate is only logarithmic unless the group is arithmetic.

\begin{figure}
\begin{tabular}{ccc}
\begin{overpic}[scale=.6]{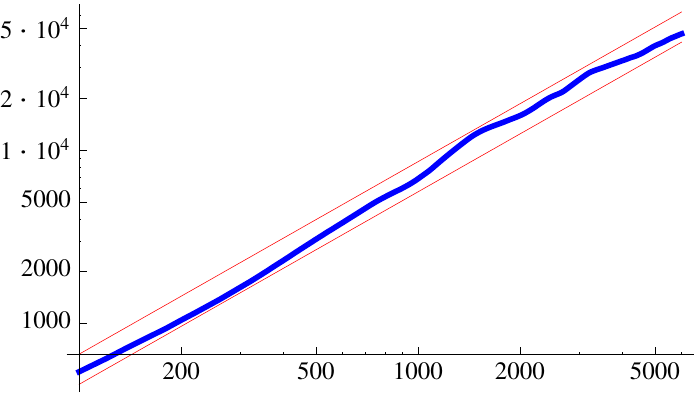}
\put(20,40){$\scriptstyle\delta\approx .1084$}
\put(20,47){$\scriptstyle X(12,12.8,13.6)$}
\end{overpic}&
\begin{overpic}[scale=.6]{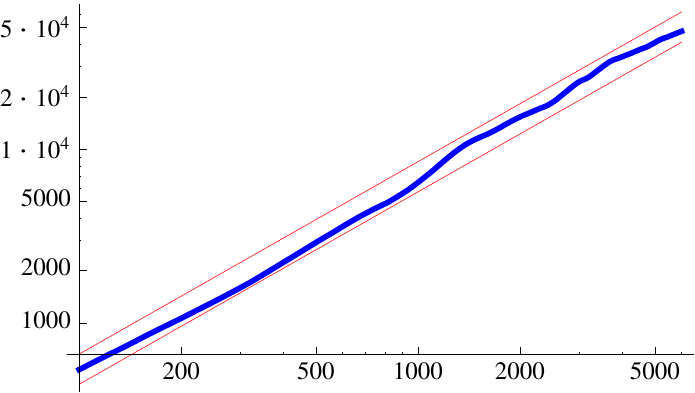}
\put(20,40){$\scriptstyle\delta\approx .1068$}
\put(20,47){$\scriptstyle X(12,13,14)$}
\end{overpic}&
\begin{overpic}[scale=.6]{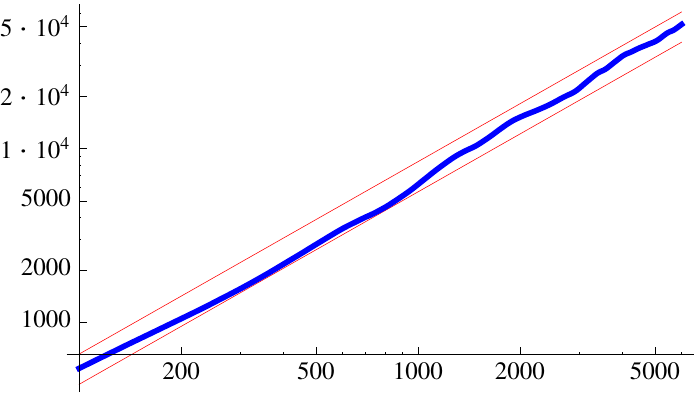}
\put(20,40){$\scriptstyle\delta\approx .1053$}
\put(20,47){$\scriptstyle X(12,13.2,14.4)$}
\end{overpic}
\end{tabular}
\caption{Log-log plots of $N_X(0,\delta;t)$ for several 3-funnel surfaces.}\label{FrWeyl12.fig}
\end{figure}
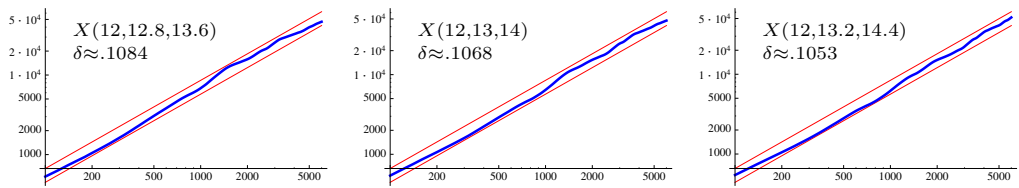

For the 3-funnel surfaces, large scale fluctuations in the growth rate seem to be a persistent feature of the vertical-strip counting function.  Figure~\ref{FrWeyl12.fig} shows a comparison across several surfaces with nearby parameters, with the bounds as conjectured in the fractal Weyl law included for comparison.   (The straight lines are graphs of $ct^{1+\delta}$, for two values of $c$ chosen arbitrarily.)  
 The results are certainly not inconsistent with the fractal Weyl law prediction, but within the range of these plots, at least, the curves are too 
unstable for any meaningful estimate of the growth rate.

\begin{figure}
\begin{center}
\begin{overpic}[scale=.9]{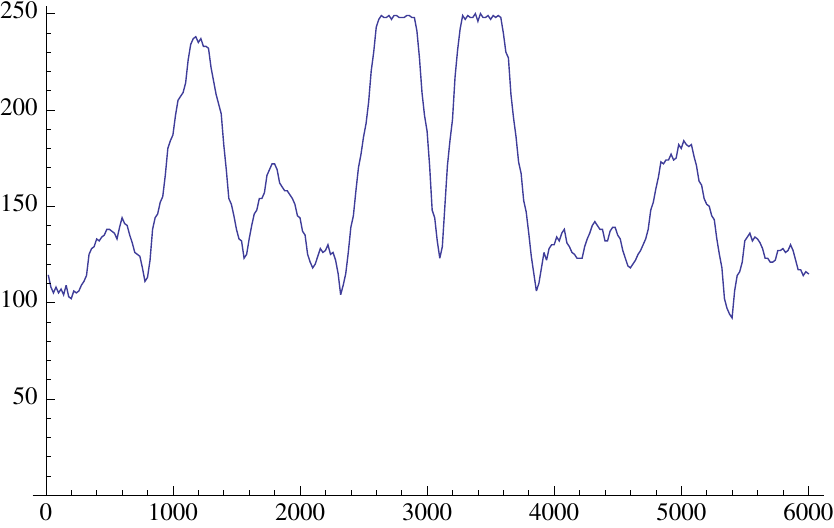}
\put(100,6){$t$}
\end{overpic}
\end{center}
\caption{The resonance count $\#\{\zeta \in \calR_X:\> \re \zeta >0,\> \abs{\im\zeta - t}\le 10\}$, for the surface $X(12,13,14)$.}\label{Nbox121314.fig}
\end{figure}

\begin{figure}
\begin{center}
\begin{overpic}[scale=.7]{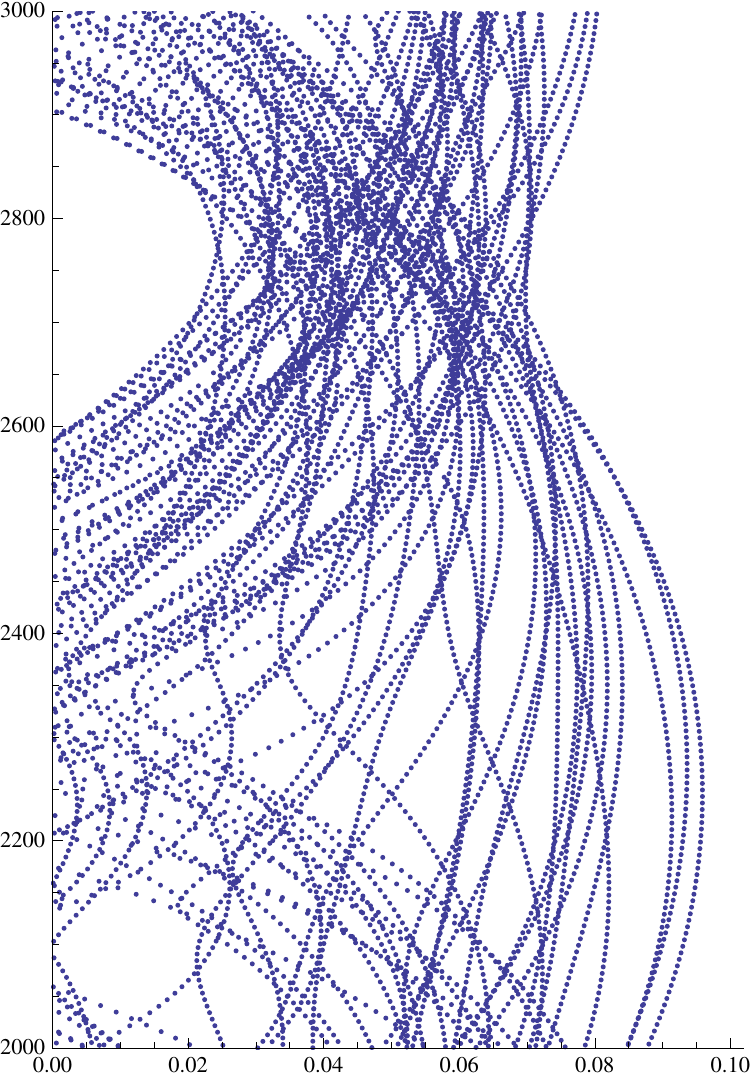}
\end{overpic}
\end{center}
\caption{Resonances with $2000 \le \im \zeta \le 3000$, for the surface $X(12,13,14)$.  The plateaus seen in Figure~\ref{Nbox121314.fig} are explained by a concentration of resonances curves, as seen between $\im s = 2600$ and $2900$.}\label{ResPlot121314.fig}
\end{figure}

We can also compare our resonance counts to the localized upper bound \eqref{glz.bound}.  Figure \ref{Nbox121314.fig} shows the plot of this resonance count in boxes of increasing imaginary height, for the surface $X(12,13,14)$ whose cumulative count was shown in the center graph in Figure~\ref{FrWeyl12.fig}.  The plateaus occurring on either side of $t=3000$ are puzzling.  To understand what is causing them, we return to a resonance plot for this surface, in Figure~\ref{ResPlot121314.fig}.  As we noted in \S\ref{resplot.sec}, on each curve of resonances the spacing is remarkably consistent.  And since the curves are nearly vertical, this means the space between imaginary parts is essentially constant on each curve.  So the fluctuations in the resonance count shown in Figure~\ref{Nbox121314.fig} are accounted for by resonance curves crossing the imaginary axis, rather than by a change in the density within individual curves.  The plateaus occur in zones where the resonance lines bunch together, creating an interval with no crossings.

\begin{figure}
\begin{center}
\begin{overpic}{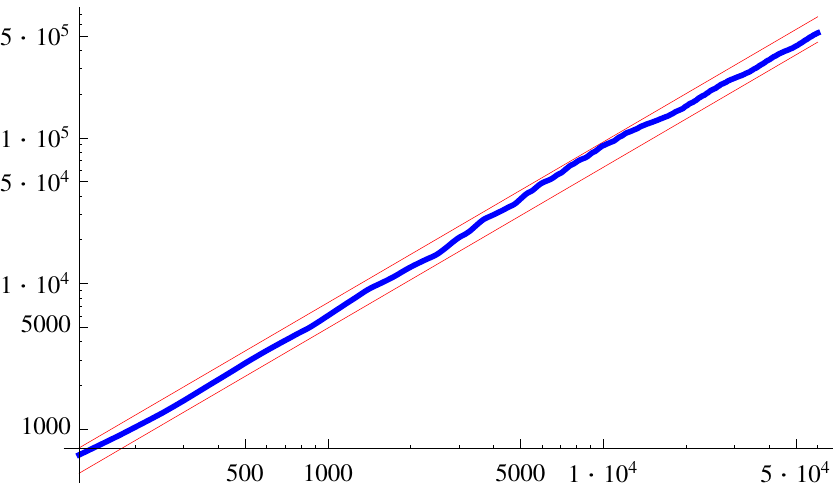}
\end{overpic}
\end{center}
\caption{Log-log plot of $N_X(0,\delta;t)$ for $X(12,14,15)$.  The outer lines show the slope predicted by the fractal Weyl law.}\label{Weyl121415.fig}
\end{figure}

\begin{figure}
\begin{center}
\begin{overpic}[scale=.8]{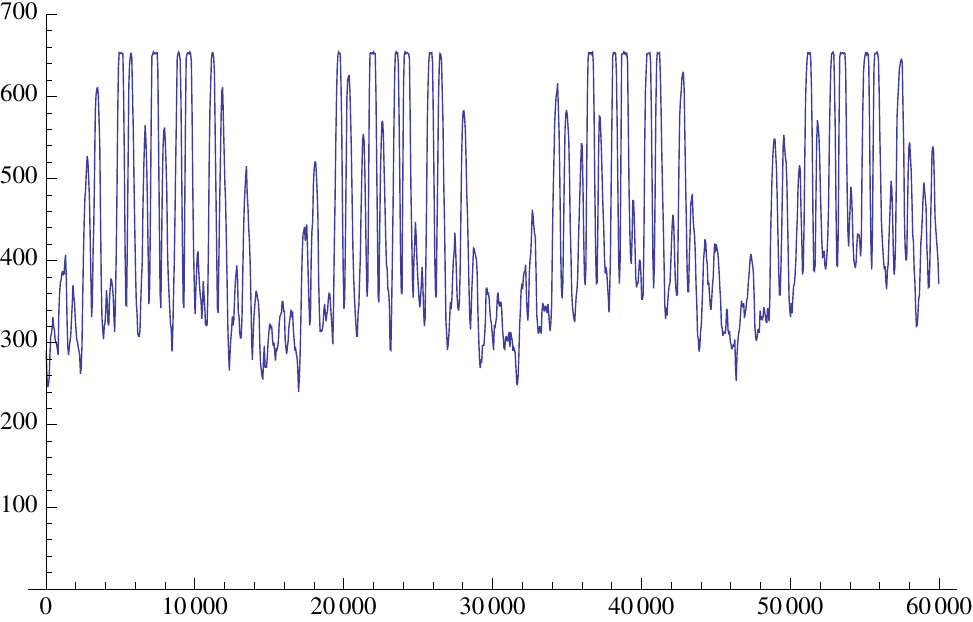}
\put(99,4){$t$}
\end{overpic}
\end{center}
\caption{The resonance count $\#\{\zeta \in \calR_X:\> \re \zeta >0,\> \abs{\im\zeta - t}\le 25\}$, for the surface $X(12,14,15)$.}\label{Nbox121415.fig}
\end{figure}

For longer values of the length parameters, we can calculate the strip counting function for much larger values of the imaginary parts.
Figure~\ref{Weyl121415.fig} shows the plot of $N_X(0,\delta;t)$ for the 3-funnel surface $X(12,14,15)$, up to $t =  60000$.
We see similar large scale oscillations even on this much greater scale.   
The plateaus in the resonance count noted above persist for large values of $t$, as shown in Figure~\ref{Nbox121415.fig}, and within the range of this plot the heights of the plateau do not vary.  Another striking feature of this graph is the quasi-periodic structure; the pattern of peaks and valleys recurs with a period somewhere around 15000.  From one cycle to the next we see a clear rise in the valleys of the distribution.  If the fractal Weyl law does hold, then it would seem to require that those valleys must eventually rise above the (apparently) fixed height of the plateaus.  Unfortunately, direct evidence of this is beyond the range of our calculations.

\begin{figure}
\begin{center}
\begin{overpic}[scale=.7]{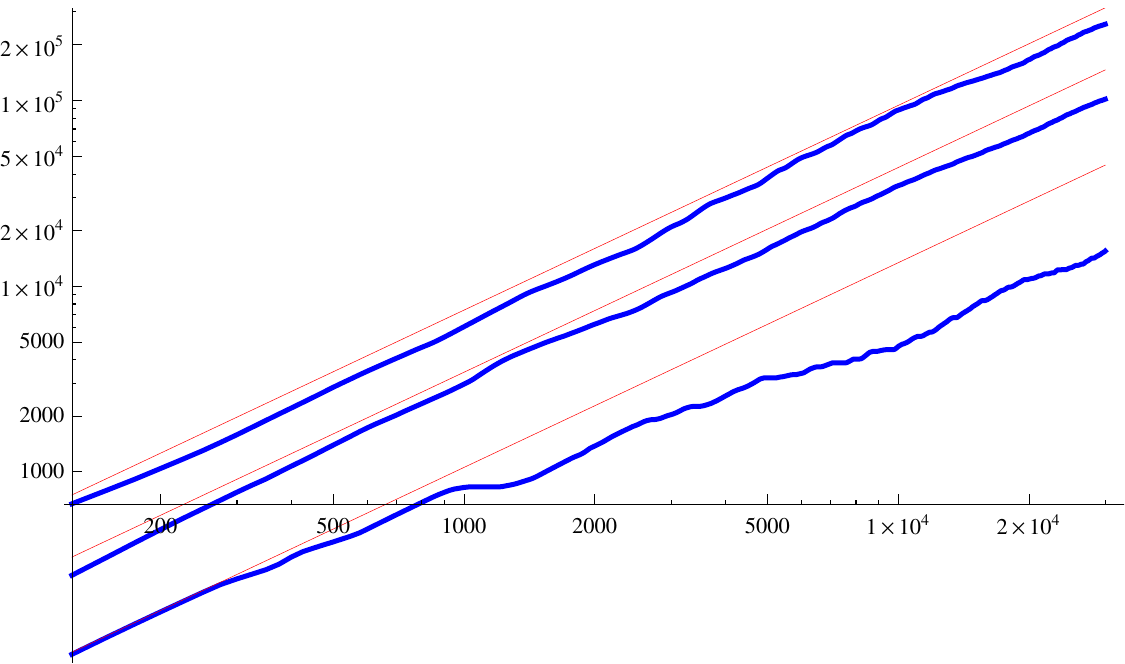}
\put(100,36){$\sigma = 3\delta/4$}
\put(100,49){$\sigma = \delta/2$}
\put(100,56.5){$\sigma = 0$}
\end{overpic}
\end{center}
\caption{For $X(12,14,15)$, a comparison of $N_X(\sigma,\delta;t)$ for several values of sigma.  The thin lines show the slope predicted by the fractal Weyl law for $\sigma = 0$.}\label{FWsigma.fig}
\end{figure}

Recently, Naud \cite{Naud:2012}, has proved, for convex co-compact $\Gamma$, that for $\sigma\ge \delta/2$ there exists a function $\tau(\sigma)$ such that  
\begin{equation}\label{nstrip.tau}
N_X(\sigma,\delta;t) = O(t^{1+\tau(\sigma)}),
\end{equation}
with $\tau(\delta/2) = \delta$ but $\tau(\sigma) < \delta$ for $\sigma > \delta/2$.  In other words, the optimal growth rate conjectured in \eqref{weyl.conj} certainly does not hold for $\sigma > \delta/2$.  We can see the evidence of this result by plotting the counting function $N_X(\sigma,\delta;t)$ for several values of $\sigma$, as in Figure~\ref{FWsigma.fig}.  For $\sigma = \delta/2$ the difference in slope is barely evident, but least-squares fit gives an exponent $1.047$, versus $1.092$ for $N_0(0,\delta;t)$.  At $\sigma = 3\delta/4$ the difference in slope is much more dramatic;  the exponent from the curve fit is only $0.888$ in this case.

\begin{figure}
\begin{tabular}{cc}
\begin{overpic}[scale=.75]{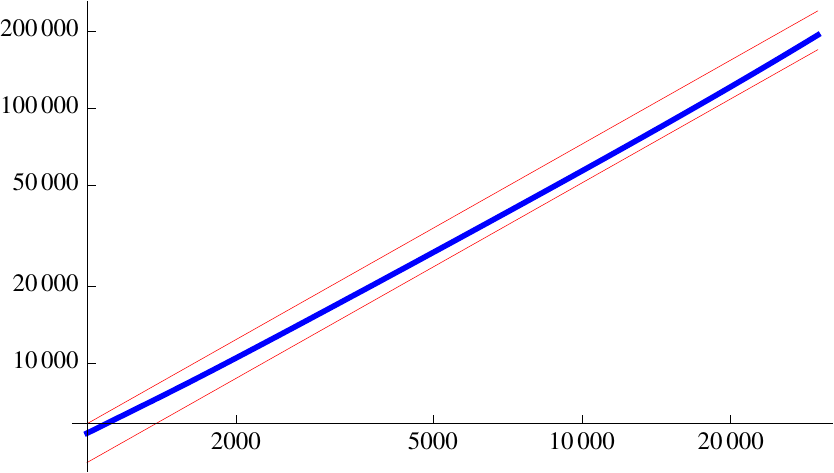}
\put(20,50){$\scriptstyle Y(12,12,\pi/2)$}
\put(20,44){$\scriptstyle \delta \approx 0.0953$}
\end{overpic} &
\begin{overpic}[scale=.75]{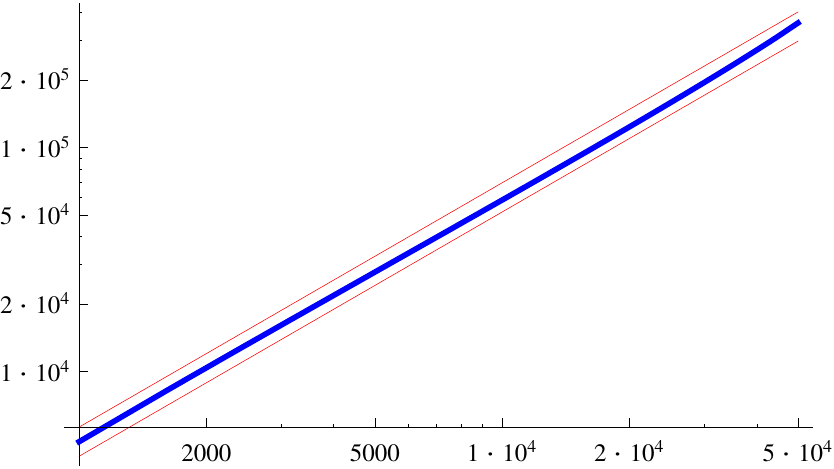}
\put(20,50){$\scriptstyle Y(12,13,\pi/2)$}
\put(20,44){$\scriptstyle \delta \approx 0.0913$}
\end{overpic}
\end{tabular}
\caption{Log-log plot of $N_X(0,\delta;t)$ for two funneled tori examples.  The thin lines show the slopes predicted by the fractal Weyl law.}\label{FrWeylTor.fig}
\end{figure}

\begin{figure}
\begin{center}
\begin{overpic}[scale=.8]{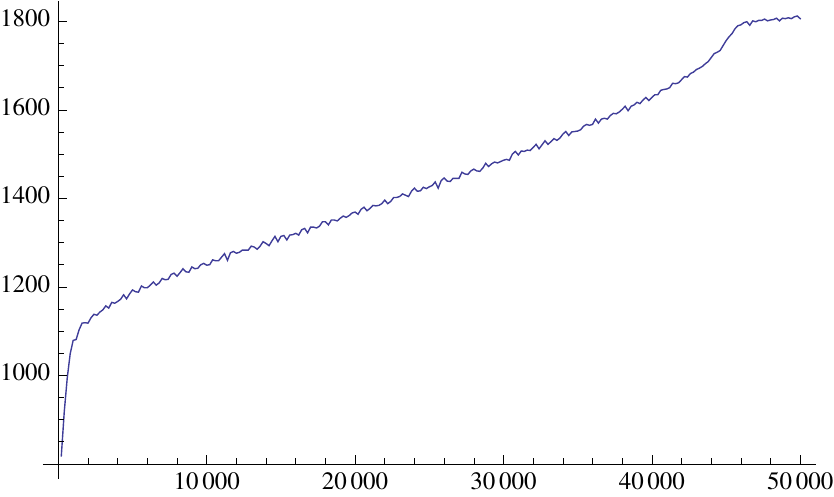}
\put(98,6){$t$}
\end{overpic}
\end{center}
\caption{The resonance count $\#\{\zeta \in \calR_X:\> \re \zeta >0,\> \abs{\im\zeta - t}\le 100\}$, for the surface $Y(12,13,\pi/2)$.}\label{Nbox1213.fig}
\end{figure}

For funneled tori the log-log plot of the resonance counting function in strips looks more uniform but still has significant long-range fluctuation in all of the examples checked.    Figure~\ref{FrWeylTor.fig} shows some examples.  Note the apparent concavity in the log-log plots.  This behavior would contradict \eqref{glz.bound} if it were to continue.   To get some idea what is happening, we examine the corresponding localized resonance count in Figure~\ref{Nbox1213.fig}.
Up to $\im s \approx 45000$, the resonance count indeed seems to have faster-than-linear growth, in disagreement with the proven $O(t^\delta)$ bound \eqref{glz.bound}.  The final plateau shows (fortunately!) that this behavior does not continue.  The resulting curve is irregular and difficult to draw any long-range conclusions from, but at least it plausibly satisfies the $O(t^\delta)$ bound.

\begin{figure}
\begin{center}
\begin{overpic}[scale=.5]{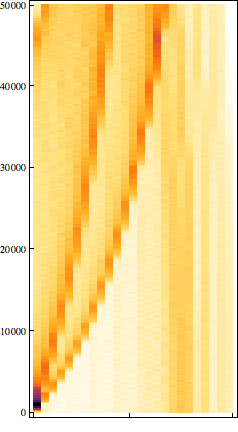}
\put(7,0){$\scriptstyle 0$}
\put(28,0){$\scriptstyle \delta/2$}
\put(53,0){$\scriptstyle \delta$}
\end{overpic}
\end{center}
\caption{Density plot of the resonance distribution for the surface $Y(12,13,\pi/2)$.}\label{Dens1213.fig}
\end{figure}

Another way to try to understand the nature of this distribution is to plot resonance densities as in Figure~\ref{Dens1213.fig}.  This plot reveals lines of higher density that evolve slowly as $\im s$ increases.  The asymptotics of the counting function would seem to be dominated by the evolution of these structures.

\section{Spectral gap}\label{spgap.sec}

In physical terms, the resonances with the greatest real part correspond to the most stable states, and hence play a dominant role in the wave asymptotics.  These wave asymptotics are related to asymptotics of the length spectrum by the trace formula, so the most stable resonances are also the most significant in the related problems of lattice point counting and asymptotics of the length spectrum.   See, for example, Naud \cite{Naud:2005b}.  In a similar way, understanding the rightmost edge of the resonance distribution plays a crucial role in various counting problems in number theory.  For examples in the context of this paper, see the work of Bourgain-Gamburd-Sarnak \cite{BGS:2011} or Bourgain-Kontorovich \cite{BK:2010}.

\begin{figure}
\begin{center}
\begin{overpic}[scale=.8]{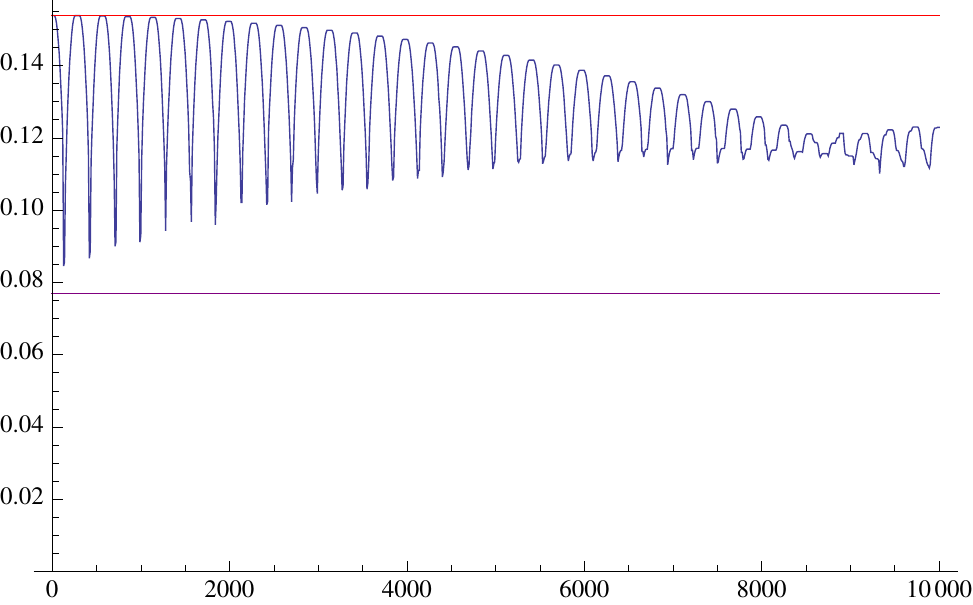}
\end{overpic}
\end{center}
\caption{Resonance envelope function $h_w(t)$ for $X(9,9,9)$ with $w = 25$.  The horizontal lines mark $\delta$ and $\delta/2$.}\label{Gap999.fig}
\end{figure}

By the positivity of $\Delta_X$, any resonances in $\re s > \tfrac12$ lie in the interval $(\tfrac12,1)$.  Hence there must always be a gap between the first resonance set and the bulk of $\calR_X$ in this case.  
Naud \cite{Naud:2005a} has shown that such a gap exists even for $\delta\le \tfrac12$:  for some $\vep>0$,
\begin{equation}\label{naud.gap}
\calR_X \cap \{\re s \ge \delta - \vep\} = \{\delta\}.
\end{equation}
One can see in various plots above that there can be resonances with $\re s$ extremely close to $\delta$.  For example, in the resonance plot for $X(12,14,15)$ in Figure~\ref{ResPlot121415.fig}, the rightmost line of resonances in nearly vertical near $\delta \approx 0.101821$;  the next resonance in the line above that point is $\zeta \approx 0.101816 + 6.28i$.  Within the class $X(\ell,\ell,\ell)$, the initial curve of resonances starting from $\delta$ becomes more vertical as $\ell$ increases.  Compare, for example, the slopes of these curves in Figures~\ref{ResPlot101.fig} and \ref{ResPlot121212.fig}.
Naud's constant $\vep$ evidently takes on arbitrarily small values, just within the class of 3-funnel surfaces.

Another version of the spectral gap problem is to find the ``essential spectral gap'', defined as
\[
G(X) := \inf \Bigr\{\sigma<\delta:\> \calR_X \cap \{\re s\ge \sigma\} \text{ is finite}\Bigr\}.
\]
Jakobson-Naud \cite{JN:2012} formulated this notion and proved that $G(X) \ge \frac{\delta(1-2\delta)}{2}$ for any convex co-compact group $\Gamma$.
They also made the following:
\begin{conj}[Essential spectral gap]\label{sp.gap}
For a convex co-compact group,
\[
G(X) = \frac{\delta}2.
\]
\end{conj}
\noindent
Note that this is consistent with the situation for finite-area surfaces, for which $\delta =1$.  In that case $G(X) \le \tfrac12$ because of the location of the critical line, and Selberg's result \cite{Selberg:1990} on accumulation of resonances shows that in fact $G(X) = \tfrac12$.  

To investigate the high-frequency behavior of the spectral gap, we define a resonance ``envelope'' function, for a given window size $w$, by
\[
h_w(t) := \max \bigl\{\re\zeta:\> \zeta \in \calR_X, \> \abs{\im \zeta - t} \le w \bigr\}.
\]
The essential gap conjecture can then be rephrased as
\[
\limsup_{t\to \infty} h_w(t) = \frac{\delta}2.
\]

Figure~\ref{Gap999.fig} shows a plot of this function for the surface $X(9,9,9)$.  The oscillatory pattern is a clear
extension of the structure seen on the left in Figure~\ref{ResPlot101.fig}.  All of the symmetric $X(\ell,\ell,\ell)$ cases share this same basic pattern for the short range oscillations.  We also contraction of these oscillations away from the line $\re s = \delta$, but this may a temporary or periodic effect.  
In the density plot for $X(12,12,12)$ shown below in Figure~\ref{Dens121212.fig}, we can see that the envelope function plot for that surface would show a very regular oscillation out to $t = 20000$, with only a very small decrease in amplitude.

\begin{figure}
\begin{center}
\begin{overpic}[scale=.8]{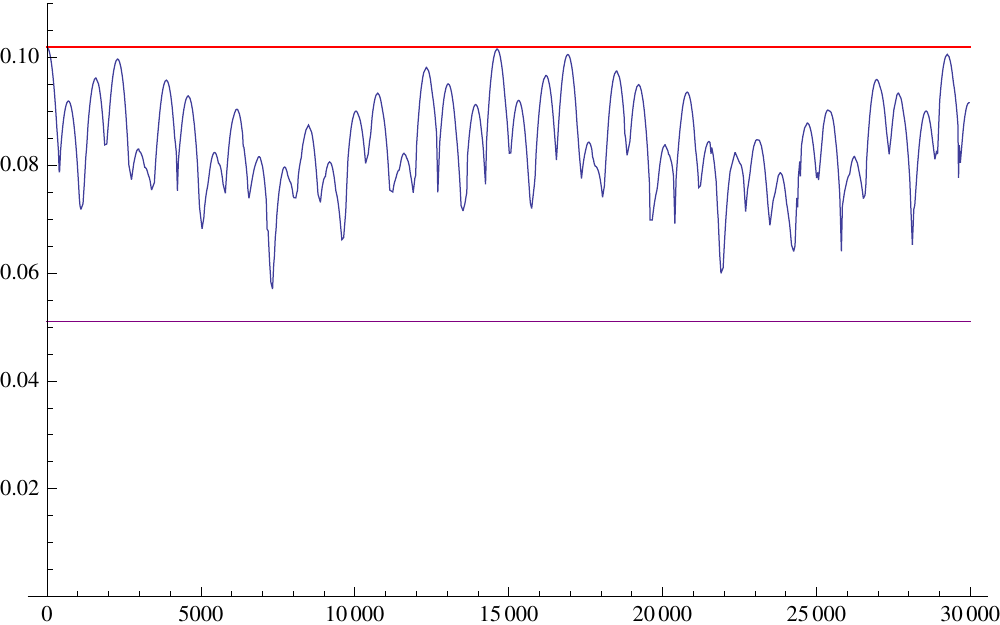}
\end{overpic}
\end{center}
\caption{Resonance envelope function $h_w(t)$ for $X(12,14,15)$ with $w = 25$.  The horizontal lines mark $\delta$ and $\delta/2$.}\label{Gap121415.fig}
\end{figure}

Figure~\ref{Gap121415.fig} shows a plot of the envelope function for $X(12,14,15)$.  The oscillatory pattern is rather striking; the shorter oscillations have a period of around 800, but there is also a clear repeat of the pattern of peaks and valleys on a larger scale, with a period of around 15000.  This of course corresponds to the periodic structure noted in Figure~\ref{Nbox121415.fig}.
The peaks of this longer oscillation occur at $\zeta = \delta \approx 0.1018$, $\zeta \approx 0.1015 + 14621.0i$, and $\zeta \approx 0.1005 + 29254.5i$.

\begin{figure}
\begin{center}
\begin{overpic}[scale=.9]{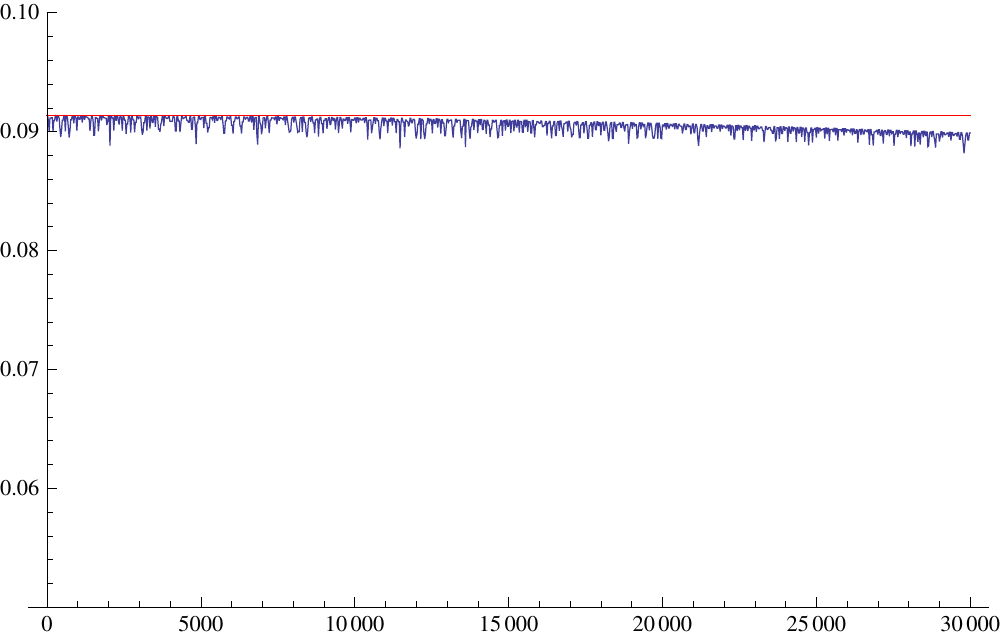}
\end{overpic}
\end{center}
\caption{Resonance envelope function $h_w(t)$ for $Y(12,13,\pi/2)$ with $w = 25$.  The horizontal line marks $\delta$.}\label{Gap1213.fig}
\end{figure}

For funneled tori, the resonances are generally more uniformly distributed, and the oscillations in the envelope function are of a much smaller amplitude.   Figure~\ref{Gap1213.fig} shows an sample plot, for the surface $Y(12,13,\pi/2)$.  
A clear opening of the gap visible on the right.   Over the course of the plot, which covers $0 \le \im s \le 30000$, the maximum real part decreases  from $\delta \approx 0.09134$ to $0.08988$.  

\section{Concentration of decay rates}\label{decay.sec}

Another interesting feature of the resonance distribution for $n$-disk scattering systems noted in Lu-Sridhar-Zworski \cite{LSZ:2003} was a concentration of the average decay rate of resonant states at half of the classical escape rate.  This observation was made in numerical calculations of the resonances for the 3-disk scattering system.  Recent experimental results for microwave scattering \cite{BWPSKZ:2012} have corroborated the existence of such a concentration.  The physical justification for this phenomenon is actually quite straightforward.  Semiclassically (i.e. for large $\abs{\im s}$), we expect the probability distributions of resonance states to exhibit, in some average sense, features of the classical flow.  In particular,  the decay rate of quantum states should be related to the classical process of escape to infinity.  Since the probability distribution is the modulus square of the wave function, the expected quantum decay rate is actually half of the classical escape rate.   
 
\begin{figure}
\begin{center}
\begin{overpic}[scale=.5]{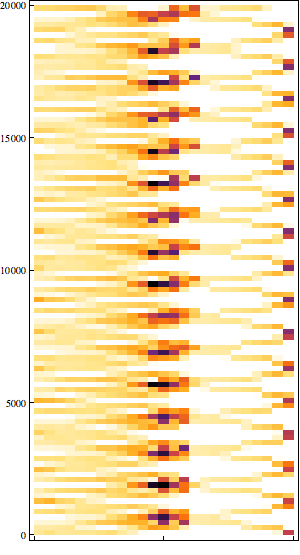}
\put(5.5,0){$\scriptstyle 0$}
\put(27.5,0){$\scriptstyle \delta/2$}
\put(52,0){$\scriptstyle \delta$}
\end{overpic}
\end{center}
\caption{Resonance density plot for $X(12,12,12)$, showing the concentration of decay rates.}\label{Dens121212.fig}
\end{figure}

To define the classical escape rate, for a small fixed $\vep>0$ and $T\ge 0$, let $U^\vep_T$ be the set of points $w \in S^*X$ for which there exists some $v \in S^*X$ in the trapped set such that $d(\phi_t(w), \phi_t(v)) <\vep$
for all $t \in [0,T]$ (distance in the Sasaki metric on $S^*X$).  In other words, $U^\vep_T$ contains those momenta for which trajectories remain close to  trapped trajectories up to time $T$.  The escape rate is defined as the decay exponent for the volume of this set,
\[
\alpha(X) := - \limsup_{T\to\infty} \frac{\log \vol(U^\vep_T)}{T} .
\]
Using the fact that the escape rate can be computed in terms of the topological pressure of the flow, Naud \cite{Naud:2005c} showed that
for a conformally compact hyperbolic surface,
\[
\alpha(X) = 1 - \delta.
\]
On the other hand, the decay rate of a resonant state at $s = \zeta$ is $\tfrac12 - \re \zeta$.  
Hence we should expect a concentration of the real parts of resonances near the line $\re s = \delta/2$.  
Naud's estimate \eqref{nstrip.tau} gives some theoretical support to this conjecture, showing at least that a reduced density of resonances for $\re s > \delta/2$. Note also that the predicted concentration is consistent with the picture for finite-area surfaces.  In that case the classical escape rate is zero, and Dyatlov \cite[eq.~(1.17)]{Dyatlov:2012} shows that under generic conditions most resonances lie very close to $\re s = \tfrac12$.

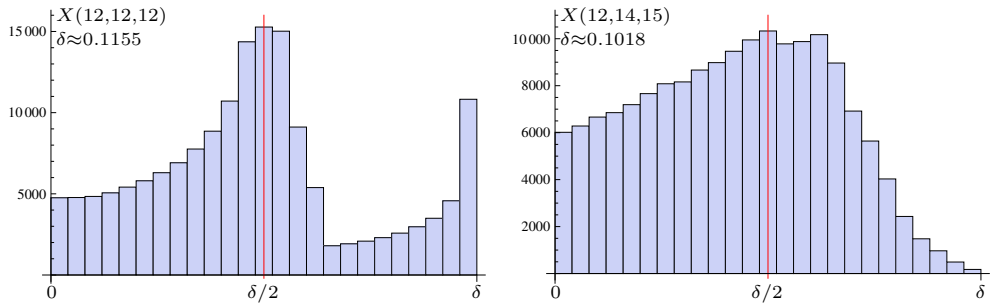
\begin{figure}
\begin{tabular}{cc}
\begin{overpic}[scale=.6]{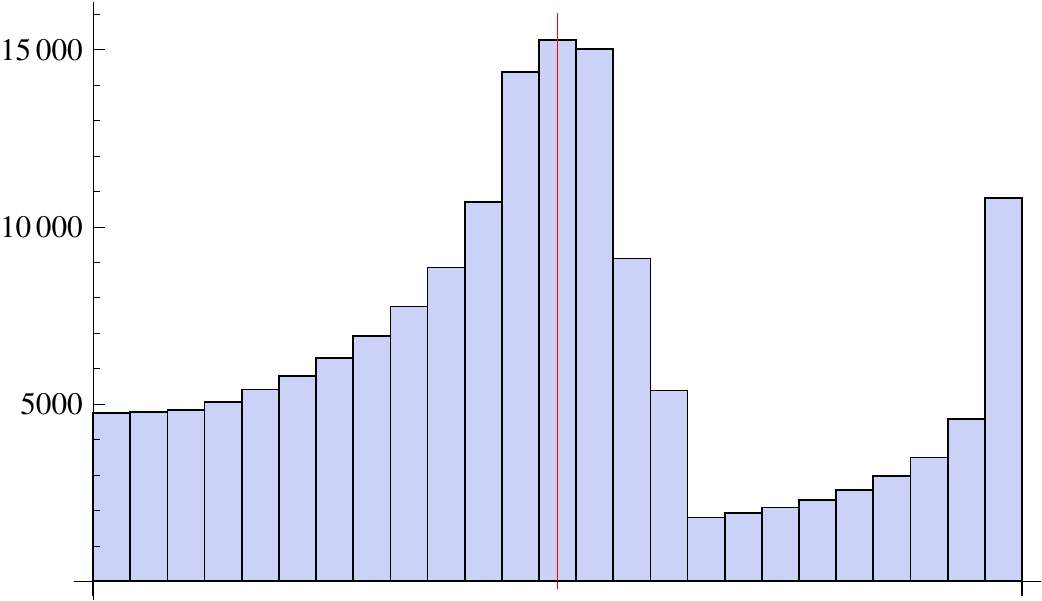}
\put(10,55){$\scriptstyle X(12,12,12)$}
\put(10,50){$\scriptstyle \delta \approx 0.1155$}
\put(8,-3){$\scriptstyle 0$}
\put(50,-3){$\scriptstyle \delta/2$}
\put(97,-3){$\scriptstyle \delta$}
\end{overpic} &
\begin{overpic}[scale=.6]{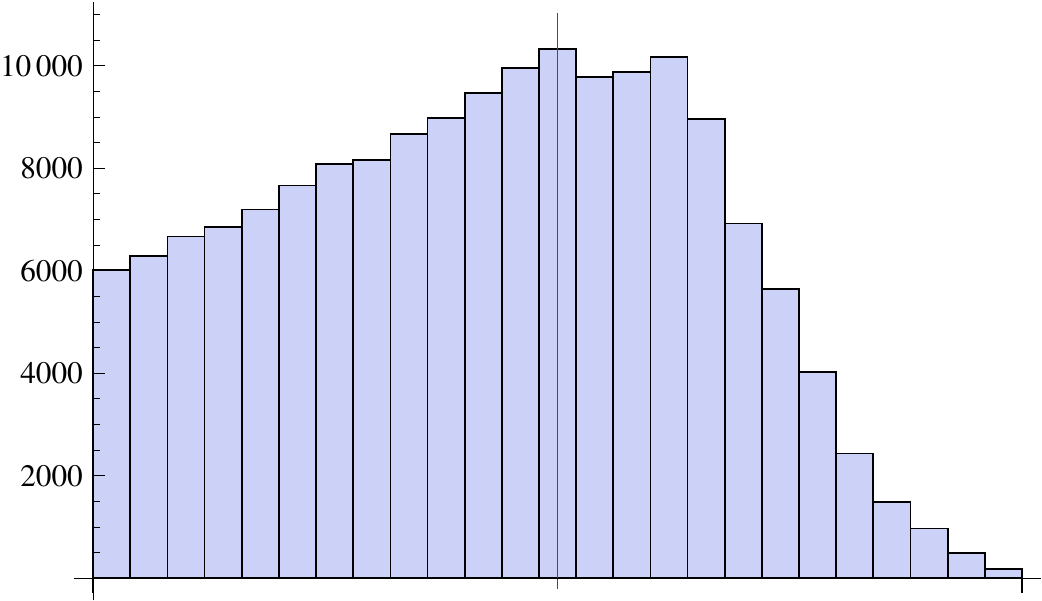}
\put(10,55){$\scriptstyle X(12,14,15)$}
\put(10,50){$\scriptstyle \delta \approx 0.1018$}
\put(8,-3){$\scriptstyle 0$}
\put(50,-3){$\scriptstyle \delta/2$}
\put(97,-3){$\scriptstyle \delta$}
\end{overpic}
\end{tabular}
\caption{Histograms of real parts of resonances with $0 \le \im s \le 20000$.  The vertical line shows $\re s= \delta/2$.}\label{RHist.fig}
\end{figure}

We have already noted the hubs of resonance curves near $\re s = \delta/2$ in Figures~\ref{ResPlot101.fig} and \ref{ResPlot121212.fig}, which shows a concentration of decay rates similar to that seen in the examples cited above.  This concentration continues over a large scale, as illustrated by the density plot in Figure~\ref{Dens121212.fig}.

\begin{figure}
\begin{center}
\begin{overpic}[scale=.55]{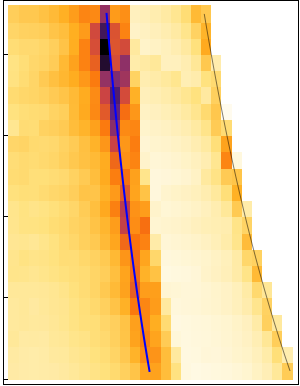}
\put(-19.5,1){$\scriptstyle X(10,10,10)$}
\put(-19.5,22){$\scriptstyle X(11,11,11)$}
\put(-19.5,43){$\scriptstyle X(12,12,12)$}
\put(-19.5,64){$\scriptstyle X(13,13,13)$}
\put(-19.5,85){$\scriptstyle X(14,14,14)$}
\end{overpic}
\end{center}
\caption{Density of real parts of resonances with $0 \le \im s \le 2000$, for a range of 3-funnel surfaces.  The lines indicate the values of $\delta/2$ and $\delta$.}\label{RDplot.fig}
\end{figure}
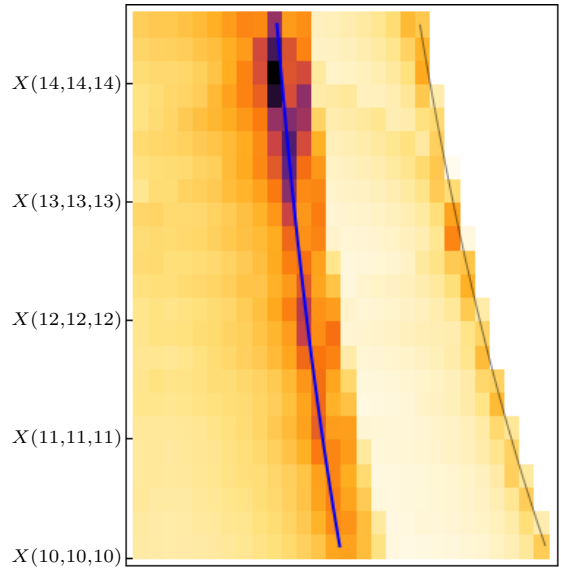

The concentration is most striking in the symmetric cases, but all of the 3-funnel surfaces seem to share the phenomenon of bunching of resonance curves near $\re s = \delta/2$ to some extent.  Figure~\ref{RHist.fig} shows a comparison of two histograms of the distribution of real parts for a portion of the resonance set.  We can track the parameter dependence of this behavior with a density plot of the real parts take over a range of length parameters, as shown in Figure~\ref{RDplot.fig}.

\begin{figure}
\begin{center}
\begin{overpic}[scale=.6]{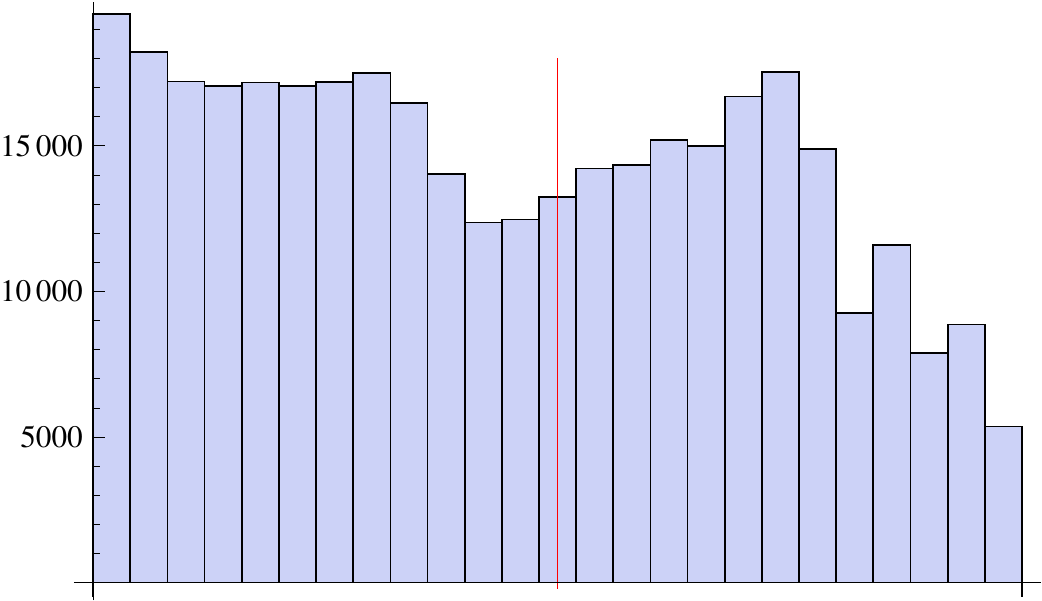}
\put(8,-3){$\scriptstyle 0$}
\put(50,-3){$\scriptstyle \delta/2$}
\put(97,-3){$\scriptstyle \delta$}
\end{overpic}
\end{center}
\caption{Histogram of the real parts of resonances with $0 \le \im s \le 30000$, for the funneled torus $Y(12,13,\pi/2)$.}\label{Rhist1213.fig}
\end{figure}

For funneled tori, the situation is quite different.  The line $\re s = \delta/2$ does not seem to be singled out in the plots of Figure~\ref{ResPlotTor.fig} or
Figure~\ref{ResPlot1012.fig}.  Figure~\ref{Rhist1213.fig} shows the histogram of real parts for a range of the resonance set for $Y(12,13,\pi/2)$, which shows no sign of the expected concentration.   In the density plot for this surface in Figure~\ref{Dens1213.fig}, 
we see that the curves on which the resonances concentrate are not coalescing near $\re s = \delta/2$, at least not within the range of the plot.  It is certainly conceivable that these curves do eventually converge on $\re s = \delta/2$, but within the range of our numerical techniques we are unable to find direct evidence of this.

\end{document}